\documentclass{article}
\usepackage{color}
\usepackage{latexsym}
\usepackage{amsmath}
\catcode`\@=11 \@addtoreset{equation}{section}

\catcode`\@=12

\newtheorem{Theorem}{Theorem}[section]
\newtheorem{Proposition}{Proposition}[section]
\newtheorem{Lemma}{Lemma}[section]
\newtheorem{Corollary}{Corollary}[section]
\newtheorem{Definition}{Definition}[section]
\newtheorem{Remark}{Remark}[section]
\newtheorem{Observation}{Observation}[section]
\newcommand{\bTheorem}[1]{
\begin{Theorem} \label{T#1} }
\newcommand{\eT}{\end{Theorem}}
\newenvironment{proof}{\textit{Proof. }}{\hfill$\Box$}

\newcommand{\bProposition}[1]{
\begin{Proposition} \label{P#1}}
\newcommand{\eP}{\end{Proposition}}

\newcommand{\bObservation}[1]{
\begin{Observation} \label{P#1}}
\newcommand{\eO}{\end{Observation}}

\newcommand{\bLemma}[1]{
\begin{Lemma} \label{L#1} }
\newcommand{\eL}{\end{Lemma}}

\newcommand{\bCorollary}[1]{
\begin{Corollary} \label{C#1} }
\newcommand{\eC}{\end{Corollary}}

\newcommand{\bDefinition}[1]{
\begin{Definition} \label{D#1} }
\newcommand{\eD}{\end{Definition}}

\newcommand{\bRemark}[1]{
\begin{Remark} \label{P#1}}
\newcommand{\eR}{\end{Remark}}

\newcommand{\bFormula}[1]{
\begin{equation} \label{#1}}
\newcommand{\eF}{\end{equation}}

\newcommand{\Ov}[1]{\overline{#1}}
\DeclareMathOperator{\suppess}{ess\, sup}
\newcommand{\DC}{C^\infty_c}

\newcommand{\vr}{\varrho}
\newcommand{\vre}{\vr_\ep}
\newcommand{\vte}{\vt_\ep}
\newcommand{\vue}{\vu_\ep}
\newcommand{\vt}{\vartheta}
\newcommand{\vT}{\varsigma}

\newcommand{\vu}{\vc{u}}
\newcommand{\vW}{\vc{W}}
\newcommand{\vn}{\vc{n}}

\newcommand{\vc}[1]{{\vec #1}}

\newcommand{\Divz}{{\rm div}_z}

\newcommand{\Divh}{{\rm div}_h}
\newcommand{\Dive}{{\rm div}_{\epsilon}}
\newcommand{\Grad}{\nabla_x}
\newcommand{\Gradz}{\nabla_z}
\newcommand{\Gradh}{\nabla_h}
\newcommand{\Grade}{\nabla_{\epsilon}}

\newcommand{\tn}[1]{\mbox {\F #1}}
\newcommand{\dx}{{\rm d} {x}}
\newcommand{\dy}{{\rm d} {y}}
\newcommand{\dt}{{\rm d} t }

\newcommand{\intO}[1]{\int_{\Omega} #1 \ \dx}

\newcommand{\ep}{\epsilon}
\newcommand{\bFrame}[1]{\bigskip
             \noindent
                         \begin{tabular}[c] {|p{.96\textwidth}|}  \hline
             \\ #1 \\  \hline \end{tabular}
             \bigskip }
\font\F=msbm10 scaled 1000
\newcommand{\R}{\mbox{\F R}}
\newcommand{\RR}{\mbox{\FF R}}
\newcommand{\RRR}{\mbox{\FFF R}}

\font\FF=msbm10 scaled 800
\font\FFF=msbm10 scaled 700
\date{}
\makeindex
\begin{document}

\title{The rotating Navier--Stokes--Fourier--Poisson system on thin domains}
\author{Bernard Ducomet$^1$,  Matteo Caggio$^2$, \v S\' arka Ne\v casov\' a$^2$,  and Milan Pokorn\' y$^3$}
\maketitle

\bigskip

\centerline{ $^1$ CEA, DAM, DIF, F-91297 Arpajon, France}
\centerline{e-mail: {\tt bernard.ducomet@cea.cz}}
\bigskip

\centerline{ $^2$ Institute of Mathematics of the Academy of Sciences of the Czech Republic}
\centerline{\v Zitn\' a 25, 115 67 Praha 1, Czech Republic}
\centerline{e-mails: {\tt caggio@math.cas.cz}, {\tt matus@math.cas.cz}}


\bigskip

\centerline{$ ^{3}$ Charles University in Prague, Faculty of Mathematics and Physics}
\centerline{Mathematical Inst. of Charles University, Sokolovsk\' a 83, 186 75 Prague 8, Czech Republic}
\centerline{e-mail: {\tt pokorny@karlin.mff.cuni.cz}}
\vskip0.25cm
\begin{abstract}
We consider the compressible Navier--Stokes--Fourier--Poisson system describing the motion of a viscous heat conducting rotating fluid
confined to a straight layer
$ \Omega _{\epsilon} = \omega \times (0,\epsilon) $, where $\omega$ is a 2-D domain. The aim of this paper is to show that the weak solutions
 in the 3D domain converge to the strong solution of the 2-D Navier--Stokes--Fourier--Poisson system on $\omega$ as $\epsilon \to 0$ on the  time interval, where the strong solution exists. We consider two different regimes in dependence on the asymptotic behaviour of the Froude number.
\end{abstract}



\medskip


\section{Introduction}
\label{i}

Our motivation for this work is a rigorous derivation of the equations describing astrophysical objects called ``accretion disk" which are thin structures
observed in various places in the universe, as e.g. Saturn's rings, protoplanetary disks, or disks found in close binary stars in which one star
 captures matter lost by its companion through a stellar wind \cite{MA} \cite{O} \cite{P}.

In fact accretion disks are among the most important objects in astronomy as they accompany the stellar formation
 and they are the main sites for planetary formation.
Namely physicists  realized already in the sixties that in a binary system (composed of two massive objects)
 mass transfer from a neutron star to its companion takes the form of an  accretion disk.

In this context, the hydrodynamical disk structure \cite{C}, \cite{Sh} comes from a subtle competition between rotation, gravitational attraction and viscosity, and
 self-gravitation seems also to play an important role in the disk evolution \cite{Pi}.

 Even though, strictly speaking, these disks are indeed three-dimensional, their size in the ``third" dimension is
 usually very small and therefore they are often modelled as two-dimensional. Our aim is to study mathematically
 rigorously the limit when the thickness tends to zero, hence the object is in the limit described as two-dimensional.

We start from the compressible Navier--Stokes--Fourier--Poisson system describing the motion of a viscous heat conducting  rotating fluid confined
to a straight layer
$ \Omega_{\epsilon} = \omega \times (0,\epsilon) $, where $\omega$ is a 2-D domain.

The motion of the fluid is described by standard fluid mechanics equations giving the
evolution of the mass density $\sigma = \sigma (t,z)$, the velocity
field $\vW = \vW(t,z)$, and the temperature $\varsigma =
\varsigma(t,z)$ as functions of the time $t$ and the spatial coordinate $z \in \Omega_{\epsilon} \subset \R^3$.
 We also assume that the flow is influenced by the gravitational force due to the fluid itself.

More precisely, the system of equations to be studied in $(0,T) \times \Omega _{\epsilon}$ reads as follows:
\bFormula{i1}
\partial_t \sigma + \Divz (\sigma \vW) = 0,
\eF
\bFormula{i2}
\partial_t (\sigma \vW) + \Divz (\sigma \vW \otimes \vW) + \sigma (\vec{\chi} \times \vW) +\Gradz p
= \Divz \tn{S} + \sigma \Gradz \Phi +\sigma \Gradz |\vec \chi \times \vec z|^2,
\eF
\[
\partial_t \left( \sigma \left( \frac{1}{2} |\vW|^2 + e \right) \right)
 + \Divz
\left(
\sigma \left( \frac{1}{2} |\vW|^2 + e \right) \vW
+  p \vW + \vc{q} - \tn{S} \vW \right)
\]
\bFormula{i3}
 =  \sigma \Gradz \Phi \cdot \vW+\sigma \Gradz |\vec \chi \times \vec z|^2 \cdot \vW,
\eF
where the body forces are represented by the gravitational force $\sigma \Gradz \Phi$,
 and the centrifugal force $\sigma \Gradz |\vec \chi \times \vec z|^2$. We namely suppose that the system is globally rotating at uniform velocity $\chi$ around the vertical direction $\vec e_3$
 and we denote $\vec \chi=\chi\vec e_3$. Therefore (see e.g. \cite{C}) the Coriolis acceleration $\sigma \vec \chi \times \vu$ and the centrifugal
force term $\sigma \Gradz |\vec \chi \times \vec z|^2$, must be included into the momentum equation.

The potential $\Phi $ obeys Poisson's equation,
\bFormula{i4}
-\Delta \Phi = 4\pi G (\alpha\sigma+(1-\alpha) g ) \ \ \  \ \mbox{in} \ (0,T) \times \Omega _{\epsilon},
\eF
where $G$ is the Newton constant and $\alpha$ is a positive parameter (see below). 
The first contribution in the right-hand side corresponds to self-gravitation while in the second one $g$ is a given function, modelling the external gravitational effect of a
central massive object (or, possibly, also of all other objects which are important to consider) attracting the system.
 Physically this may correspond to the gravitational attraction by a black hole (in the classical limit) or by a central star in the case of circumstellar disks \cite{Sh}.

Here and hereafter, we assume that the function $\sigma$ is extended by zero outside of $\Omega_\ep$.
Supposing further that  $g$ is such that the integral below converges, we have
\[
\Phi(t,z)=G\int_{\RR^3}K(z-y)\big(\alpha\sigma(t,y)+(1-\alpha)g(y)\big)\ \dy,
\]
where $K(z)=\frac{-1}{|z|}$, and the parameter $\alpha$ may take the values 0 or 1: for $\alpha=1$ the self-gravitation is present and for
$\alpha=0$ the gravitation only acts as an external field (astrophysicists often consider self-gravitation of accretion disks as small compared to the external
 attraction described by $g$, see e.g. \cite{P}). We will explain the choice of the gravitational potential later. Since we also work with  $\Gradz \Phi$, we have further to
  assume that
\[
\int_{\RR^3}|\nabla K(z-y)| (\alpha \sigma (t,y) + (1-\alpha)g(y)|)\ \dy <\infty.
\]

We consider  the no slip boundary condition holds on the boundary $\partial\omega \times (0,\epsilon)$
(on the lateral part of the domain)
\bFormula{p1}
\vW |_{\partial \omega \times (0,\epsilon)}=\vc{0},
\eF
 and the slip boundary condition on the other parts of the boundary  $\omega \times \{0,\epsilon\}$ (the top and bottom part of the layer)
\bFormula{p1bis}
 \vW \cdot \vn|_{ \omega \times \{0,\epsilon\} } = 0,\ \ [\tn{S}\vec n ] \times \vec n|_{  \omega \times \{0,\epsilon\}}=\vec 0.
\eF
Note that  we have on $\omega \times \{0,\epsilon\}$ $\vec n = \pm \vec e_3$, hence the first condition in (\ref{p1bis}) can be rewritten as
\bFormula{p1bisa}
W_{3} = 0 \mbox { on } \omega \times \{0,\epsilon\}.
\eF

Let us emphasize that we imposed a slip condition on the boundary $\omega \times \{0,\epsilon\}$ in order to avoid difficulties in
 passing to the ``infinitely thin" disk limit.

We suppose finally that the no heat flux condition holds on $\partial\Omega$
\bFormula{i11}
 \vc{q} \cdot \vc{n} |_{\partial \Omega_\ep} = 0.
\eF

Equations are completed with the initial conditions
\bFormula{in1}
\sigma(0,z) =  \widetilde {\sigma}_{0,\epsilon }(z),\ \ \vc W (0,z) = \widetilde  {\vc W}_{0,\epsilon},
\ \ \vT (0,z) =\widetilde  {\vT}_{0,\epsilon }(z)\ \ \mbox{for}\ z \in \Omega _{\epsilon}.
\eF

The thermodynamic pressure $p$ and the specific internal energy $e$ are given functions of the density $\sigma$ and the temperature $\vT$. They are  interrelated through Maxwell's relation
\bFormula{i5}
\frac{\partial e(\sigma,\vT)}{\partial \sigma} = \frac{1}{\sigma ^2} \left( p(\sigma,\vT) - \vT \frac{\partial p(\sigma,\vT)}{\partial \vT} \right).
\eF
Furthermore, $\tn{S}$ is the viscous part of the stress tensor determined by
\bFormula{i6}
\tn{S}=\tn{S}(\vT,\Gradz \vW) = \mu \left( \Gradz \vW + (\Gradz \vW)^T - \frac{2}{3} \Divz \vW \right) + \eta \ \Divz \vW \ \tn{I},
\eF
where the shear viscosity coefficient $\mu = \mu(\vT) > 0$ and the bulk viscosity coefficient $\eta = \eta(\vT) \geq 0$ are effective functions of the
 temperature, and the upper index $T$ denotes the transposed matrix.

Similarly, the heat flux $\vc q$ is given by Fourier's law
\bFormula{i7}
\vc q = \vc q(\vT,\Gradz \vT) = - \kappa(\vT) \Gradz \vT,
\eF
with the heat conductivity coefficient $\kappa(\vT) > 0$.

Using standard tools and methods from the mathematical compressible fluid dynamics (see \cite{FEINOV}) it is possible to establish the existence of a weak solution to our system (\ref{i1}--\ref{i7}) under physically reasonable assumptions on the form of the pressure $p(\sigma,\vT)$ and the functions $\mu(\vT)$, $\eta(\vT)$ and $\kappa(\vT)$. Our aim is to verify that under reasonable assumptions on the data of the problem,  for $\epsilon \to 0^+$, the solution tends to the strong solution of the corresponding two-dimensional system (presented below).
The essential tool in this paper will be the relative entropy inequality as developed by E. Feireisl, A. Novotn\' y and co-workers
in \cite{FeBuNo}, \cite{FN} and \cite{FNS}
(for utilisation of relative entropy inequalities in other contexts we refer to works of
C. Dafermos \cite{D}, P. Germain \cite{G}, A. Mellet and A. Vasseur \cite{MV} and L. Saint-Raymond \cite{SR}).

The aim of this paper is to extend recent results of E. Feireisl, A. Novotn\' y and co-workers (see \cite {BFN}, \cite{MN}) to the rotating Navier--Stokes--Fourier
 system in presence of the gravitational forces.

The paper is organized as follows.
In Section \ref{m}, we collect our main assumptions on the form of the pressure $p(\sigma,\vT)$, the viscosity coefficients and the heat conductivity  and introduce
the concept of weak solution to our original problem. We also introduce the two-dimensional target system and present the main result of the paper which claims that the weak solution to the original system converge for $\epsilon \to 0^+$ to the strong solution of the target system (on the life span of this solution). Section \ref{nsl} contains the derivation of the relative entropy inequality and  Section \ref{pr} the proof of the convergence result.

\section{Hypotheses and the target system}
\label{m}

\subsection{Thermodynamic concept}
\label{mI}

As mentioned above, we use the temperature and the density as the main thermodynamic quantities and assume that the thermodynamic potentials are given functions of them.
Hypotheses imposed on constitutive relations and transport
coefficients are motivated by the general existence theory for the
Navier--Stokes--Fourier system developed in \cite[Chapter 3]{FEINOV}. As we cannot handle the monoatomic gas model when the self-gravitation is included, we use the generalization of the model as considered e.g. in \cite{NP}; see also \cite{KrNePo_ZAMP}.

We consider the pressure and the internal energy in the form
\bFormula{m0a}
p(\sigma,\vT) = p_1(\sigma,\vT) + \frac a3 \vT^4,
\eF
\bFormula{m0b}
e(\sigma,\vT) = e_1(\sigma,\vT) + a \frac {\vT^4}{\sigma},
\eF
where
\bFormula{m0c}
p_1(\sigma,\vT) = (\gamma-1) \sigma e(\sigma,\vT)
\eF
with $\gamma >1$. Note that $\gamma=\frac 53$ corresponds exactly to the model of the monoatomic gas, as considered in \cite{FEINOV}.
Recalling the Maxwell relation (\ref{i5}), it follows under some regularity assumptions on the functions $p_1$ and $e_1$ that
\bFormula{m1}
p_1(\sigma, \vT) = \vT^{\frac{\gamma}{\gamma-1}} P \left( \frac{\sigma}{\vT^{\frac{1}{\gamma-1}}} \right),
\eF
where $P: [0, \infty) \to [0, \infty)$ is a given function with the following properties:
\bFormula{m2}
P \in C^1([0, \infty))\cap C^2((0,\infty)), \ P(0) = 0, \ P'(Z) > 0 \ \mbox{for all} \ Z \geq 0,
\eF
\bFormula{m3}
0 < \frac{ \gamma P(Z) - P'(Z) Z }{Z} \leq c <\infty  \ \mbox{for all} \ Z > 0,
\eF
\bFormula{m4}
\lim_{Z \to \infty} \frac{P(Z)}{Z^{\gamma}} = p_\infty > 0.
\eF

The component $\frac{a}{3} \varsigma^4$ represents the effect of
``equilibrium'' radiation pressure (see \cite{DFN} for the motivations).

Using the Maxwell relation (\ref{i5}), the corresponding non-radiative part of the
internal energy $e$ is
\bFormula{m5}
e_1(\sigma, \vT) = \frac{1}{\gamma-1} \ \frac{\vT^{\frac{\gamma}{\gamma-1}}}{\sigma}\ P
\left( \frac{\sigma}{\vT^{\frac{1}{\gamma-1}}} \right).
\eF
Note also that the Maxwell relation is equivalent with the Gibbs relation
\bFormula{m4a}
\frac{1}{\vartheta} \Big(De(\sigma,\vT) + p(\sigma,\vT)D\Big(\frac{1}{\sigma}\Big)\Big) = Ds(\sigma,\vT)
\eF
which defines a new quantity, the specific entropy. Due to the form of the pressure and the internal energy the specific entropy is given by
\bFormula{m6}
s(\sigma, \vT) = s_1(\sigma,\vT) + \frac{4a}{3} \frac{\vT^3}{\sigma},
\eF
with
\bFormula{m6a}
s_1(\sigma,\vT) = M \left( \frac{\sigma}{\vT^{\frac{1}{\gamma-1}}} \right), \quad M'(Z) = - \frac{1}{\gamma-1} \frac{ \gamma P(Z) - P'(Z) Z } {Z^2} < 0, \quad \lim_{Z\to +\infty} M(Z) = 0.
\eF
Note that it is possible to show that
\bFormula{Ent1}
s_1(\sigma,\varsigma) \leq C(1+ |\ln \sigma|)
\eF
in the set $\sigma \in (0,\infty)$, $\varsigma \in (0,1)$, and
\bFormula{Ent2}
s_1(\sigma,\varsigma) \leq C(1+ |\ln \sigma| + \ln \varsigma)
\eF
in the set $\sigma \in (0,\infty)$ and $\varsigma \in (1,\infty)$.

The transport coefficients $\mu$, $\eta$, and $\kappa$ are continuously differentiable functions of the absolute temperature such that
\bFormula{m7a}
0 < c_1 (1 + \vT) \leq \mu (\vT), \ \mu'(\vT) < c_2, \ 0 \leq \eta(\vT) \leq c_3(1 + \vT),
\eF
\bFormula{m8a}
0 < c_4 (1 + \vT^3) \leq \kappa(\vT) \leq c_5 (1 + \vT^3)
\eF
for any $\vT >0 $.
For the sake of simplicity, we consider the particular case
\bFormula{m7}
\mu (\vT)= \mu _0 + \mu_1\vT , \quad  \mu_0 , \mu_1 >0, \qquad \eta \equiv 0,
\eF
and
\bFormula{m8}
\kappa(\vT)=\kappa _0 + \kappa_2\vT ^2 + \kappa _3\vT ^3, \quad \kappa _i>0, i=0,2,3.
\eF
As a matter of fact, the total energy balance (\ref{i3}) is not
suitable for the weak formulation since, at least according to the
recent state-of-art, the term $\tn{S} \vW$ is not controlled on the
(hypothetical) vacuum zones of vanishing density. Following
\cite{FEINOV}, we replace (\ref{i3}) by the internal energy balance (it is possible at least formally, for smooth solutions, by subtracting the balance of the kinetic energy, i.e. the momentum equation tested by $\vW$)
\bFormula{m15}
\partial_t (\sigma e) + \Divz (\sigma e \vW) + \Divz \vc q = \tn{S}(\vT,\Gradz \vW): \Gradz \vW - p \Divz \vW.
\eF
Furthermore, dividing (\ref{m15}) by $\vT$ and using Gibbs' relation (\ref{m4a}), we may rewrite (\ref{m15}) as the entropy equation
\bFormula{m16}
\partial_t (\vr s) + \Divz (\sigma s \vW) + \Divz \left( \frac{\vc{q}}{\vT} \right) =
\frac{1}{\vT} \Big( \tn{S}(\vT,\Gradz \vW) : \Gradz \vW - \frac{ \vc{q} \cdot \Gradz \vT }{\vT} \Big).
\eF

\subsection{Formal scaling analysis}
\label{mII}
We introduce the change of variables
\bFormula{c1}
(z_h, z_3)\in \Omega _{\epsilon} \mapsto (x_h, x_3) \in \Omega_1=\Omega, \mbox { where } x_h = (x_1,x_2), z_h=(z_1,z_2).
\eF
Note that for $z \in \Omega_\ep$ we have $x \in \Omega_1 =: \Omega$.
For $z$ and $x$ interrelated by \eqref{c1} we introduce new unknowns defined as follows: the density $\vr(t,x) = \sigma(t,z)$, the velocity $\vu(t,x)= \vc{W}(t,z)$ and the thermodynamic temperature $\vt(t,x) = \vT(t,z)$. Similarly we set $\phi(t,x) = \Phi(t,z)$.

In order to identify the appropriate limit regime, we perform a general scaling, denoting by
 $L_{ref},
\ell_{ref},
\ T_{ref},
\ U_{ref}$,
$\rho_{ref},
\ p_{ref},
\ \mu_{ref},$
the reference hydrodynamical quantities (horizontal and vertical lengths, time, velocity, density, pressure, viscosity) and by
$\chi_{ref}$ the reference rotation velocity.

We denote by $Sr:=\frac{L_{ref}}{T_{ref}U_{ref}}, \
Ma:=\frac{U_{ref}}{\sqrt{\rho_{ref}p_{ref}}}, \
Re:=\frac{U_{ref}\rho_{ref}L_{ref}}{\mu_{ref}},$ $
Fr:=\frac{U_{ref}}{\sqrt{G\rho_{ref}L_{ref}}},$
$Ro:=\frac{U_{ref}}{\chi_{ref}L_{ref}},$ $Pe:=\frac{U_{ref}p_{ref}L_{ref}}{\vt_{ref}\kappa_{ref}}$, the Strouhal, Mach,
Reynolds, Froude, Rossby and P\'eclet numbers. We set
\[
\epsilon:=\frac{\ell_{ref}}{L_{ref}}.
\]

We also denote
\[
\nabla _{\epsilon}= (\nabla _h, \frac{1}{\epsilon}\partial _{x_3}), \quad \nabla _h = (\partial _{x_1},\partial _{x_2}),
\]
\[ \Dive  \vc u = \Divh \vu_h + \frac{1}{\epsilon}\partial _{x_3}u_3,\ \ \vc u_h = (u_1,u_2),\ \ \ \Divh \vc u_h =
 \partial _{x_1}u_1+ \partial _{x_2}u_2,
\]
\[
\Delta_{\ep} = \partial^2_{x_1 x_1} + \partial^2_{x_2 x_2} + \frac{1}{\epsilon^2}\partial^2_{x_3 x_3}.
\]
Note also that
\[
\Grade \phi(t,x)
=\ep G\int_{\Omega}\alpha \vr(t,\xi)\frac{(x_1-\xi_1,x_2-\xi_2, \ep(x_3-\xi_3))}{\big(|x_h-\xi_h|^2 + \epsilon^2 |x_3-\xi_3|^2\big)^{\frac 32}}\ {\rm d}\xi  \\
\]
\bFormula{A0}
+ G\int_{\RR^3}(1-\alpha)g(y)\frac{(x_1-y_1,x_2-y_2, \ep x_3-y_3)}{\big(|x_h-y_h|^2 + (\epsilon x_3-y_3)^2\big)^{\frac 32}}\ {\rm d}y =: \ep \alpha \vc{E}_1 + (1-\alpha) \vc{E}_2 =: \vc{E},
\eF
and
\bFormula{A1}
\Grade |\vec\chi \times \vec x|^2 = (\nabla_h |\vec \chi \times \vec x|^2,0) = \frac{\chi}{\sqrt{x_1^2 + x_2^2}} (x_1,x_2,0).
\eF

Performing the scaling
the continuity equation reads now
\bFormula{i1bis}
 Sr\ \partial_t \vr + \Dive (\vr \vu)= 0,
\eF
the momentum equation
\[
 Sr\ \partial_t  (\vr \vu)
 + \Dive (\vr \vu \otimes \vu)
+ \frac{1}{Ma^2}\ \Grade p(\vr,\vt)
 + \frac{1}{Ro}\vr (\vec \chi \times \vu)
\]
\bFormula{i2bis} =  \frac{1}{Re}\ \Dive \tn{S}(\Grade\vu)
+\frac{1}{Fr^2}\vr \vc{E}
 +\frac{1}{Ro}\vr \nabla_{\epsilon} |\vec \chi \times \vec x|^2,
\eF
and the entropy balance
\bFormula{i3bis}
Sr \partial_t (\vr s )
 + \Dive \left(\vr s \vu\right)
+ \frac{1}{Pe}\ \Dive \left( \frac{\vc{q(\vt,\Grade \vt)}}{\vt}\right)
= \Sigma,
\eF
with
\[
\Sigma=
\frac{1}{\vt} \left( \frac{Ma^2}{Re} \tn{S}(\vt,\Grade\vu):\Grade \vu -\frac{1}{Pe}\frac{\vc q(\vt,\Grade \vt)\cdot \Grade\vt}{\vt}\right).
\]

We take $Fr =\ep^\beta$ and set all other numbers equal to $1$. We consider in fact only two (the most interesting limit) cases: $\beta = \frac 12$  and $\beta =0$. If $\beta = \frac 12$, we take $\alpha =1$, i.e. we assume only the self-gravitation and neglect the gravity due to the other objects, and if $\beta =0$, we take $\alpha =0$, i.e. we neglect the self-gravitation effects of the fluid and consider only the gravitational force due to the other objects. Indeed, if $\beta =0$, we could also include into the primitive system the self-gravitation of the fluid, however, in the limit $\ep \to 0^+$ the corresponding part of the gravitational potential would disappear. We comment on this case below Theorem \ref{main}.

We obtain the new (rescaled) system of equations
\bFormula {c2}
\partial_t \vr + \Dive  (\vr \vu) = 0,
\eF
\bFormula {c3}
\partial_t (\vr \vu) + \Dive (\vr \vu \otimes \vu) + \Grade p(\vr,\vt) +  \vr (\vec \chi \times \vu)
 = \Dive \tn{S}(\vt, \nabla_\epsilon  \vc u) + \ep^{-2\beta} \vr \vc{E}+\vr\Grade |\vec \chi\times \vec x|^2,
\eF
\bFormula {c4}
\partial_t (\vr s) + \Dive (\vr s \vu) + \Dive \left( \frac{\vc{q(\vt,\Grade \vt)}}{\vt} \right) =
\frac{1}{\vt} \Big( \tn{S}(\vt, \Grade \vu) : \Grade \vu - \frac{ \vc{q}(\vt,\Grade \vt) \cdot \Grade \vt }{\vt} \Big),
\eF
with initial conditions
\bFormula{c5}
\vr(0,x) = \vr_{0,\epsilon }(x), \, \vc u(0,x) = \vc u_{0,\epsilon}, \, \vt (0,x) = \vt_{0,\epsilon }(x),\ \  x \in \Omega,
\eF
where
\[
 \vr_{0,\epsilon }(x) =\widetilde \sigma_{0,\epsilon }(z),
\ \ \vc u_{0,\epsilon }(x) =\widetilde \vc{W} _{0,\epsilon }(z),
\ \ \vt_{0,\epsilon }(x) =\widetilde \vT_{0,\epsilon }(z),
\]
and where $\vc{E}$ is given by \eqref{A0}.

The rescaled problem (\ref{c2}--\ref{c5}) is completed with the boundary conditions
\bFormula{p2}
\vu |_{\partial \omega \times (0,1)}=\vc{0},
\eF
\bFormula{p2bis}
\vu\cdot \vec n|_{  \omega \times \{0,1\}}=0,\ \ \ [\tn{S}(\vt, \Grade \vu )\vec n ] \times \vec n|_{  \omega \times \{0,1\}}=0,
\eF
and
\bFormula{p2ter}
\Grade \vt \cdot \vec n|_{  \omega \times (0,1)}=0.
\eF
Recall that due to the shape of $\Omega$ we have
\bFormula{p2terb}
u _3 =0 \mbox { on } \omega \times \{0,1\}.
\eF

\subsection{Weak formulation}
\label{mIII}

To present the weak formulation of our problem, we consider the functional space
\[
 W^{1,2}_{0,n}(\Omega; \R^3 ) = \{\vu \in W^{1,2}(\Omega; \R^3 );\,  \vu \cdot \vn|_{\omega \times \{0,1\} } = 0, \vu|_{\partial \omega \times (0,1)}= \vc{0} \}.
\]

In the weak formulation of the Navier--Stokes--Fourier system the equation of continuity  (\ref{c2})  is replaced by its
(weak) renormalized version (cf. \cite{DL}) represented by the family of
integral identities
\bFormula{m12}
\int_0^T \intO{ \Big( \big( \vr + b(\vr) \big) \partial_t \varphi + \big(\vr + b(\vr) \big) \vu \cdot \Grade \varphi + \big( b(\vr) - b'(\vr) \vr \big) \Dive \vu \varphi \Big)}
\ \dt
\eF
\[
 = - \intO{ \big( \vr_{0,\ep} + b(\vr_{0,\ep}) \big) \varphi (0, \cdot) }
\]
satisfied for any $\varphi \in \DC([0, T) \times \Ov{\Omega})$,
and any $b \in C^\infty([0,\infty))$, $b' \in \DC([0, \infty))$, where
(\ref{m12}) includes as well the initial condition $\vr(0,
\cdot) = \vr_0$.

Similarly, the momentum equation (\ref{c3}) is replaced by
\[
\int_0^{T} \intO{ \big( \vr \vu \cdot \partial_t \vec\varphi + (\vr \vu \otimes \vu) : \Grade \vec\varphi
 - \vr (\vec \chi \times \vu)\cdot\vec\varphi+  p(\vr,\vt) \Dive \vec\varphi \big) } \ \dt
\]
\bFormula{m13}
-\int_0^{\tau} \intO{ \left(\tn{S}(\vt,\Grade \vu) : \Grade \vec\varphi  - \ep^{-2\beta}\vr \vc{E}\cdot \vec\varphi -\vr \Grade |\vec \chi \times \vec x|^2 \cdot\vec\varphi \right)}\ \dt
\eF
\[
=-\intO{ \vr_{0,\epsilon} \vu_{0,\epsilon}\cdot \vec\varphi (0, \cdot) }
\]
for any $\vec\varphi \in \DC([0, T) \times \overline{ \Omega}; \R^3 )$ such that $\vec \varphi |_{[0,T]\times \partial \omega \times (0,1)}=\vec{0}$ and
 $ \varphi _3 |_{[0,T]\times \omega \times \{0,1\}}=0$, where (\ref{m13}) includes the initial condition  $\vr \vu (0,
\cdot) = \vr_{0,\epsilon}\vu_{0,\epsilon}$.

The entropy equation must be replaced by the entropy inequality. Its weak formulation is the following
\bFormula{m17}
\int_0^T \intO{ \left(  \vr s(\vr,\vt) \partial_t \varphi + \vr s(\vr,\vt) \vu \cdot \Grade \varphi
 + \frac{ \vc q(\vt,\Grade \vt) }{\vt} \cdot \Grade \varphi \right) } \ \dt
\eF
\[
\leq - \intO{ \vr_{0,\ep} s(\vr_{0,\ep},\vt_{0,\ep}) \varphi (0, \cdot) }
- \int_0^T \intO{
\frac{1}{\vt} \left( \tn{S}(\vt,\Grade \vu) : \Grade \vu - \frac{\vc{q(\vt,\Grade \vt)} \cdot \Grade \vt }{\vt} \right) \varphi } \ \dt,
\]
for any $\varphi \in \DC([0,T) \times \Ov{\Omega})$, $\varphi \geq 0$.
Since we replace the entropy equation by the inequality, the system  must be supplemented with
the total energy balance
\bFormula{m18}
\intO{ \left( \frac{1}{2} \vr |\vu|^2 + \vr e(\vr, \vt) \right)( \tau, \cdot) }
\eF
\[
= \intO{ \left( \frac{1}{2\vr_{0,\ep} }|(\vr_{0,\ep} \vu_{0,\ep} |^2 + \vr_{0,\ep} e(\vr_{0,\ep},\vt_{0,\ep})   \right) }
+\int_0^{\tau}\int_{\Omega}\left[\ep^{-2\beta}\vr \vc{E}\cdot \vu +\vr \Grade |\vec \chi \times \vec x|^2 \cdot \vu\right]\ \dx \ \dt.
\]
Finally, the gravitational force is given by  (\ref{A0}).

\bDefinition{d1}
\label{d1}
 We say that $(\vr, \vu, \vt)$ is a weak solution of problem (\ref{c2}--\ref{c5}) if
\[
\vr \geq 0, \ \vt > 0 \ \mbox{for a.a.}\ (t,x) \in (0,T)\times \Omega,
\]
\[
\vr \in C_{weak}([0,T]; L^{\gamma}(\Omega)), \ \vr\vu \in  C_{weak}([0,T]; L^{\frac{2\gamma}{\gamma+1}}(\Omega;\R^3)),
\]
\[\vu \in L^2(0,T; W^{1,2}_{0,n}(\Omega; \R^3)),
\]
\[
\vt \in L^\infty(0,T; L^4(\Omega)) \cap  L^2(0,T; W^{1,2}(\Omega)),
\]
and if $\vr$, $\vu$, $\vt$ satisfy the integral identities (\ref{m12}), (\ref{m13}), (\ref{m17}),  together with the total  energy balance  (\ref{m18}), and integral representation of the gravitational force (\ref{A0}).
\eD

We denote for fixed $\epsilon>0$ by the lower index $\ep$ the solution to our problem on $\Omega$ (i.e., after rescaling). Recall that if $\alpha = 0$, we neglect the self-gravitation.

\bProposition{le1}
\label{le1}
Suppose the thermodynamic functions $p$, $e$, $s$
satisfy hypotheses (\ref{m0a}--\ref{m6a}),  the transport coefficients $\mu$, $\lambda$, $\kappa$ comply with (\ref{m7}--\ref{m8})
and the stress tensor is given by (\ref{i6}). Let $\gamma > \frac 32$ if $\alpha =0$ or $\gamma > \frac {12}{7}$ if $\alpha =1$. Let $g$ be such that $g \in L^p(\R^3)$ for $p=1$ if $\gamma >6$ and $p= \frac{6\gamma}{7\gamma-6}$ for $\frac 32<\gamma \leq 6$.

Suppose the initial data satisfy
\bFormula{l1}
\intO{ \left( \frac{1}{2} \vre |\vue|^2 + \vre e(\vre, \vte)\right)(0, \cdot) }
\eF
\[
\equiv \intO{ \left( \frac{1}{2 \vr_{0,\ep}} |\vr_{0,\ep} \vu_{0,\ep} |^2 + \vr_{0,\ep} e(\vr_{0,\ep},\vt_{0,\ep})\right) } \leq E_0,
\]
\[
\intO{ \vre s(\vre, \vte)(0, \cdot) } \equiv \intO{ \vr_{0,\ep} s(\vr_{0, \ep},\vt_{0,\ep}) } \geq S_0.
\]

Then problem (\ref{c2}--\ref{c5}) with boundary conditions (\ref{p2}--\ref{p2ter}) admits at least one weak solution in the sense of  Definition \ref{d1}.
\eP

\begin{proof}
The existence of weak solutions to the above problem can be deduced from the works of Feireisl et al. \cite{DFPS}, \cite{F},\cite{FN} and \cite{EF71}. In fact, we fix $\ep>0$, construct a weak solution in $\Omega_\ep$ and then  rescale the solution. Note that due to the form of the total energy inequality (\ref{m17}) we need stronger assumptions on $\gamma$ and $g$ than in \cite{DFPS}. However, this form will be important in Section \ref{pr}. Note that in real situations we may assume that $g$ is compactly supported, which removes the fact that the assumptions on the integrability of $g$ are in a certain sense more restrictive for large $\gamma$.
\end{proof}

\subsection{Target system}
\label{mIV}

The aim of our paper is to investigate the limit process $ \epsilon \to 0^+$ in the system of equations (\ref{c2}--\ref{c4}) under the assumptions that initial
data $(\vr_{0,\epsilon}, \vc u_{0,\epsilon}, \vt _{0,\epsilon})(x) $ converge in a certain sense to $(r_0,\vc U_0, \Theta_0)(x)$, where $\vc{U}_0 = (\vc{V}_{0},0)$ and all quantities are independent of $x_3$.

We therefore expect that the sequence $(\vr _{\epsilon}, \vc u_{\epsilon},\vt_{\epsilon})(t,x)$
of weak solutions to (\ref{c2}--\ref{c4}) converge to $(r,\vc{U}, \Theta)$, $\vec{U} = [\vc{V},0]$ has zero third component and all quantities are independent of $x_3$.

Then the triple $(r(t,x_h),{\vec V}(t,x_h), \Theta (t,x_h))$
solves the 2D compressible heat conducting rotating fluid in the domain $(0,T) \times\omega $. However, the target system depends on the gravitational potential, i.e. on the choice of $\beta$ (and $\alpha$).
\bFormula{t1}
 \partial_t r + \Divh  (r \vc{V}) = 0,
\eF
\bFormula {t2}
r\partial_t \vc{V} + r \vc{V} \cdot \Gradh \vc{V} + \Gradh p(r,\Theta) +  r (\vec \chi \times \vc{U})
 = \Divh \tn{S}_h(\Theta,\Gradh  \vc{V} ) +  r \Gradh \phi_h+ r \Gradh |\vec\chi\times \vec x|^2,
\eF
\bFormula {t3}
r \partial_t  s + r \vc{V} \cdot \Gradh s + \Divh \left( \frac{\vc{q}_h(\Theta,\Gradh \Theta)}{\Theta} \right) =
\frac{1}{\Theta} \Big( \tn{S}_h(\Theta,\nabla_h \vc{V}) : \Gradh \vc{V} - \frac{ \vc{q}_h(\Theta,\Gradh \Theta)
 \cdot \Gradh \Theta }{\Theta} \Big).
\eF
Above,
\bFormula{t4a}
\phi_h(t,x_h)=G \int_{\omega} \frac{r(t,y_h)}{|x_h-y_h|}\ \dy_h
\eF
if $\alpha =1$ and $\beta =\frac 12$, and
\bFormula{t4}
\phi_h(t,x_h)= G \int_{\RR^3} \frac{g(y)}{\sqrt{|x_h-y_h|^2 + y_3^2}} \ \dy
\eF
if $\alpha =0$ (and $\beta =0$).
Moreover, $\vc{q}_h(\Theta,\Gradh \Theta) = -\kappa(\Theta) \nabla_h \Theta$. Further, we have $\vc{V}|_{\partial \omega \times (0,T)} = \vc{0}$,  $\vc{q}_h\cdot \vc{n}|_{\partial \omega \times (0,T)} =0$, and
\bFormula {t5}
\tn{S}_h(\Theta,\Gradh \vc{V})= \mu \left(\Gradh \vc{V} + (\Gradh \vc{V})^T - \Divh \vc{V}\right) +\Big(\eta + \frac{\mu}{3}\Big)\Divh \vc{V} \, \tn{I}_h,
\eF
where $\tn{I}_h$ is the unit tensor in $\R^{2\times 2}$.

Note that since we consider only smooth solutions to (\ref{t1}--\ref{t5}), the  entropy balance is in fact expressed as equality and it is equivalent with either the internal or the total energy balance (in the sense of (\ref{i3})).

Our aim now is to prove that solutions of (\ref{c2}--\ref{c4}) converge in a certain sense to the unique strong solution
of (\ref{t1}--\ref{t5}) as explained below.

>From classical results of Matsumura and Nishida \cite{MN2} we know that the target system admits a unique global strong solution provided the initial data are close to a stationary solution.

Another possible result is the existence of local in time smooth solution, i.e. a smooth solution on possibly short time interval $(0,T_*)$, see e.g. \cite{Ta} for such a type of result. More precisely

\bProposition{pro1}
\label{pro1}
 Let $E$ be a given positive constant.
Suppose that  the functions $p \in C^2((0,\infty)^2)$, $\mu, \, \eta, \, \kappa \in C^1 (0,\infty)$ and that
\bFormula{u1}
r_0 \in W^{2,2}(\omega), \inf _{\omega } r_0 >0, \vc V _0 \in W^{3,2}(\omega;\R^2)\cap W^{1,2}_0(\omega;\R^2), \Theta _0 \in   W^{3,2}(\omega), \inf_{\omega} \Theta_0 >0.
\eF
Assume that the compatibility condition
 \bFormula{comp1}
 \frac{1}{r_0}\Big(\Gradh p(r_0,\Theta_0) +r_0 (\vec \chi \times \vec V_0) - \Divh \tn{S}_h(\Theta_0,\Gradh  \vc V_0 )
-r_0 \Gradh \phi_h-r_0 \Gradh |(\vec\chi \times \vec x)|^2\Big)|_{\partial\omega}=\vc{0},
\eF
holds.
Then: \newline
1. {\it (Local solution)} There exists a positive constant $T_*$ such that $(r,\vec V ,\Theta )$ is the unique classical solution to problem (\ref{t1}--\ref{t5}) with the boundary conditions
\bFormula{t6a}
\vc{V}|_{\partial \omega} = \vc{0},
\eF
\bFormula{t7a}
\frac{\partial\Theta}{\partial \vec{n}}|_{\partial \omega} = 0.
\eF 
and the initial conditions $(r_0, \vc{V}_0, \Theta_0)$ in  $(0,T)\times \Omega$ for any $T<T_*$ such that
\bFormula{u3}
r\in C([0,T];W^{3,2}(\omega)) \cap C^1([0,T];W^{2,2}(\omega)) ,
\eF
\bFormula{u3bis}
\vc V \in C([0,T];W^{3,2}(\omega,\R^2))\cap C^1([0,T];W^{1,2}(\omega,\R^2)),
\eF
\bFormula{u3ter}
\Theta \in C([0,T];W^{3,2}(\omega))\cap C^1 ([0,T];W^{1,2}(\omega)).
\eF
2. {\it (Global solution)} Let $r_0,\vc V_0,\Theta_0,\chi$ be such that for a sufficiently small $\varepsilon >0$
\bFormula{closetoeq}
\|r_0 - \overline{ r}, \vec V_0,\Theta _0 - \overline {\Theta} \|_{3,2} + |\chi| \leq \varepsilon,
\eF
where $(\overline {r},\vc{0},\overline {\Theta})$ is a stationary solution to (\ref{t1}--\ref{t5}) with the boundary condition
\bFormula{t7}
\frac{\partial \overline{\Theta}_0}{\partial \vec{n}}|_{\partial \omega} = 0.
\eF

Then  for any $T_*<+\infty$ there exists a global unique strong solution
to problem
 (\ref{t1}--\ref{t5}) with the boundary conditions  (\ref{t6a}--\ref{t7a}) and the initial conditions $(r_0, \vc{V}_0, \Theta_0)$
in the class (\ref{u3}--\ref{u3ter}).
\eP
\begin{proof}
It follows from \cite[Theorem 1.1]{MN2} and \cite{Ta} with slight modifications due to the rotation and the self-gravitation. Moreover, in the former case, we also use that the potential part of the external force does not need to be small.
\end{proof}
\bigskip

\subsection{Main result}
\label{mV}

We introduce the relative entropy functional
\bFormula{ref}
{\mathcal I}\left(\vr_{\varepsilon},\vc u_{\varepsilon}, \vt_{\varepsilon}| r,\vc U, \Theta \right)= \int_{\Omega}
\left(
\frac{1}{2}\ \vr_{\varepsilon} |\vc u_{\varepsilon}-\vc U|^2
+{\mathcal E}\left(\vr_{\varepsilon},\vt_{\varepsilon}|r,\Theta\right)
\right)(\tau,\cdot)\ \dx,
\eF
where for
\[
H_{\widetilde\vt}(\vr,\vt)=\vr e(\vr,\vt)-\widetilde\vt\vr s(\vr,\vt)
\]
we define
\[
{\mathcal E}(\vr,\vt|\widetilde\vr,\widetilde\vt):=
H_{\widetilde\vt}(\vr,\vt)
-\partial_{\vr}H_{\widetilde\vt}(\widetilde\vr,\widetilde\vt)(\vr-\widetilde\vr)
-H_{\widetilde\vt}(\widetilde\vr,\widetilde\vt).
\]

Our main result reads
\bTheorem{main}
\label{main}
Suppose   the thermodynamic functions $p$, $e$, $s$
satisfy hypotheses (\ref{m0a} - \ref{m6a}),  the transport coefficients $\mu$, $\lambda$, $\kappa$ comply with (\ref{m7}--\ref{m8})
and the stress tensor is given by (\ref{i6}).

Suppose that (if $\alpha =0$)
\bFormula{l0}
\int_{\RR^3} \frac{g(y) y_3}{(\sqrt{(x_h-y_h)^2 + y_3^2})^3} \ \dy = 0 \quad \forall x_h \in \omega.
\eF

Let $r_0, \vc V_0,\Theta _0$ satisfy assumptions of Proposition \ref{pro1}  and let $T_{*} >0$ be the  time interval
of existence of the strong solution to problem (\ref{t1}-\ref{t5}) corresponding to $r_0, \vc  V_0, \Theta _0$ determined in Proposition \ref{pro1}. Let either $\beta=0$, $\alpha =0$, $\gamma >\frac 32$ and $g \in L^p(\R^3)$ with $p=1$ for $\gamma >6$ and $p= \frac{6\gamma}{7\gamma-6}$ for $\gamma \in (\frac 32, 6]$, or $\beta = \frac 12$, $\alpha =1$ and $\gamma \geq \frac{12}{5}$.

Let $(\vre, \vue, \vte )$ be a sequence of weak solutions to the 3-D compressible Navier--Stokes--Fourier(--Poisson) system (\ref{m12}--\ref{m18}) with (\ref{A0}), emanating
 from the initial data  $(\vr_{0,\ep}$, $\vu_{0,\ep}$, $\vt _{0,\ep})$.

 Suppose that
\bFormula{l1a}
{\mathcal I}(\vr_{0,\epsilon},\vu_{0,\epsilon}, \vt_{0,\epsilon},|r_{0},\vc U_{0},\Theta_{0} )\to 0,
\eF
where $\vc U_0 =(\vc V_0,0)$.
Then
\bFormula{l2}
\suppess_{t \in [0,T_{max}]} {\mathcal I}(\vr_{\epsilon},\vu_{\epsilon}, \vt_{\epsilon}| ,r_,\vc U,\Theta )\to 0,
\eF
\bFormula{l2b}
\vue \to \vc{U} = (\vc{V},0) \qquad \mbox{ strongly in } L^2(0,T; W^{1,2}(\Omega;\R^3)),
\eF
\bFormula{l2c}
\vte \to \Theta  \qquad \mbox{ strongly in } L^2(0,T; W^{1,2}(\Omega)),
\eF
\bFormula{l2d}
\log \vte \to \log \Theta  \qquad \mbox{ strongly in } L^2(0,T; W^{1,2}(\Omega)),
\eF
where the triple $(r,\vc{V},\Theta)$ satisfies the 2-D compressible Navier--Stokes--Fourier(--Poisson) system (\ref{t1}--\ref{t4}) with the boundary conditions
(\ref{p2}) and (\ref{p2ter}) on the time interval
 $[0,T]$ for any $0<T<T_*$.
\eT

\bRemark{re0}
\label{re0}
Indeed, we may also include for $\beta =0$ in the primitive system  the part corresponding to the self-gravitation of the fluid. However, passing with $\ep \to 0^+$, this term tends to zero, therefore we do not consider it here, as it would lead to an additional restriction for $\gamma$.
\eR

\bRemark{re0a}
\label{re0a}
Condition (\ref{l0}) is the necessary condition for the validity of the 2-D system, as it means that the gravitational force in the $x_3$-direction in $\omega$ is zero. 
Indeed, we may replace this condition by adding another force term to the right-hand side of (\ref{t2}) which would compensate this gravitational force. The proof is the same as in our case.
\eR

\bRemark{re1}
\label{re1}
>From (\ref{l2}) it follows, in addition to Theorem \ref{main}
\bFormula{w1}
\vre \to r \ \mbox{in} \ C_{\rm weak} ([0,T]; L^{\gamma}(\Omega)),\ \ \ \vre \to r \ \mbox{a.a. in} \ (0,T) \times \Omega.
\eF
\eR
\bRemark{re0ab}
Let us mention that we assume more stronger assumptions in Theorem \ref{main} in the case $\alpha =1,\beta =\frac 12$ than in Proposition \ref{Ple1} since we need estimates independent of $\epsilon$ see (\ref{est1}).
\eR

\bCorollary{co1}
Suppose that the thermodynamic functions $p$, $e$, $s$
satisfy hypotheses (\ref{m0a}--\ref{m6a}), that the transport coefficients $\mu$, $\lambda$, $\kappa$ comply with (\ref{m7}--\ref{m8})
and the stress tensor is given by (\ref{i6}).

Assume that $(\vr_{\epsilon,0}, \vc u_{\epsilon,0},\vt_{\epsilon,0}) $, $\vr_{\epsilon,0}\geq 0$, $\vt_{\epsilon,0} \geq 0$
satisfy
\bFormula{le3} \int_0^1 \vr_{0,\epsilon}(x) \ {\rm d} x_3 \to r_0\ \ \  \mbox {weakly in } L^1(\omega),
\eF
 \bFormula{le4}
\int_0^1 \vr_{0,\epsilon}(x) \vc u_{0,\epsilon} \ {\rm d} x_3 \to  \vc V_0\ \ \   \mbox {weakly in } L^1(\omega;\R^2),
\eF
  \bFormula{le5}
\int_0^1 \vt_{0,\epsilon}(x)  \ {\rm d} x_3 \to \Theta _0\ \ \   \mbox {weakly in } L^1(\omega),
\eF
where $r_0,\vc w_0, \Theta_0$ belong to the regularity class (\ref{u1}), and
\[
\int_{\Omega}
\left(
\frac{1}{2} \vr_{0,\epsilon} |\vc u_{0,\varepsilon}|^2
+\vr_{0,\varepsilon}e(\vr_{0,\varepsilon},\vt_{0,\ep})
\right)\ \dx
\]
\bFormula{le6}
\to
\int_{\omega}\left(
\frac{1}{2} r_{0} |\vc V_{0}|^2
+r_{0}e(r_0,\Theta_0)
\right)\ \dx_h.
\eF
Let $(\vr_{\epsilon}, \vc u_{\epsilon},\vt_{\epsilon}) $  be a sequence of weak solutions to the 3-D compressible Navier--Stokes--Fourier system (\ref{c2}--\ref{c5})
 emanating from the initial data  $(\vr_{0,\epsilon}, \vc u_{0,\epsilon},\vt_{0,\epsilon}) $.

 Then (\ref{l2}--\ref{l2d}) hold.
\eC

\section{Relative entropy inequality}
\label{nsl}

We introduce the relative entropy inequality which is satisfied by any weak solution of the Navier--Stokes--Fourier--Poisson system.
We follow here the ideas from the paper \cite{FN}.

Let us consider a triple $\{r,\vc U,\Theta\}$,  smooth functions such that $r$ and $\Theta$ are bounded below away from zero
 in $[0,T]\times \Omega$, $ \vc U\Big|_{\partial\omega\times (0,1)}=\vec{0}$, and $U_3\Big|_{\omega \times \{0,1\}}=0$.

We call {\it ballistic free energy} the thermodynamic potential given by
\[
H_{\widetilde\vt}(r,\Theta)=r e(r,\Theta)-\widetilde{\vt}r s(r,\Theta),
\]
where $(r,\Theta,\widetilde{\vt})$ are sufficiently regular function so that the definition as well as the computations below are meaningful. We denote $(\vr_\ep, \vu_\ep, \vt_\ep)$ the weak solution to (\ref{m12}--\ref{m18}) with \eqref{A0} for a fixed $\ep>0$.
The {\it relative entropy} is then defined by
\[
{\mathcal E}(\vr,\vt|\widetilde\vr,\widetilde\vt):=
H_{\widetilde\vt}(\vr,\vt)
-\partial_{\vr}H_{\widetilde\vt}(\widetilde\vr,\widetilde\vt)(\vr-\widetilde\vr)
-H_{\widetilde\vt}(\widetilde\vr,\widetilde\vt).
\]

Using as a test function in (\ref{m12}) with $b \equiv 0$ the function $\varphi=\frac{1}{2}\ |\vc U|^2$ and using that $\vc{U}$ is sufficiently smooth, we arrive after standard manipulations at
\[
\int_{\Omega}\frac{1}{2}\ (\vr_{\varepsilon}|\vc U|^2)(\tau,\cdot)\ \dx
-\int_{\Omega}\frac{1}{2}\ \vr_{0,\varepsilon}|\vc U(0,\cdot)|^2\ \dx
\]
\bFormula{mass1}
=\int_0^{\tau}\int_{\Omega}
\vr_{\varepsilon}
\left(
\vc U\cdot\partial_t\vc U+\vc u_{\varepsilon}\cdot \Grade \vc U\cdot\vc U
\right)\ \dx\ \dt .
\eF
Next, using as test function in (\ref{m13}) $\phi=\vc U$, we get similarly as above
 \[
\int_{\Omega}  \vr_{\varepsilon} \vc u_{\varepsilon}(\tau,\cdot) \cdot \vc U(\tau,\cdot)\ \dx
-\int_{\Omega}  \vr_{0,\varepsilon} \vc u_{0,\varepsilon} \cdot \vc U(0,\cdot)\ \dx
\]
\bFormula{mom1}
\begin{array}{l}
\displaystyle =\int_0^{\tau}\int_{\Omega}
\left(
\vr_{\varepsilon} \vc u_{\varepsilon}\cdot \partial_t\vc U+(\vr_{\varepsilon} \vc u_{\varepsilon}\otimes \vc u_{\varepsilon}):\Grade \vc U
+p(\vr_{\varepsilon},\vt_\ep)\ \Dive \vc U-\tn{S}(\vt_{\varepsilon},\Grade \vu_\ep):\Grade \vc U
\right)\ \dx\ \dt
\\[10pt]
\displaystyle  - \int_0^{\tau}\int_{\Omega}\left( \vr_{\varepsilon}(\vc \chi \times \vc u_{\varepsilon})\cdot \vc U - \ep^{-2\beta}
\vr_{\varepsilon}\vc{E} \cdot \vc{U} - \vr_{\varepsilon}\Grade |\vc \chi \times \vc x|^2\cdot  \vc U\right)\ \dx\ \dt.
\end{array}
\eF
Combining (\ref{mass1}), (\ref{mom1}) and (\ref{m18}), we get
\[
\int_{\Omega}
\left(\frac{1}{2}\ \vr_{\varepsilon} |\vc u_{\varepsilon}-\vc U|^2
+\vr_{\varepsilon}e(\vr_{\varepsilon},\vt_\ep)\right)(\tau,\cdot)\ \dx
\]
\[
=\int_{\Omega}
\left(\frac{1}{2}\ \vr_{0,\varepsilon} |\vc u_{0,\varepsilon}- \vc U(0,\cdot)|^2
+\vr_{0,\varepsilon}e(\vr_{0,\varepsilon},\vt_{0,\ep})\right)\dx \]
\[ +\int_0^{\tau}\int_{\Omega}\left(\ep^{-2\beta} \vr_{\varepsilon}\vc{E}\cdot (\vc u _{\varepsilon}-\vc U)-
\vr_{\varepsilon}(\vc \chi \times \vc u_{\varepsilon})\cdot (\vc u_{\varepsilon}-\vc U)+ \vr_{\varepsilon}\Grade |\vc \chi \times \vc x|^2 \cdot (\vc u_{\varepsilon} - \vc U)\right)\ \dx \ \dt
\]
\bFormula{enertot1}
+\int_0^{\tau}\int_{\Omega}
\left(
\left(
\vr_{\varepsilon}  \partial_t\vc U
+\vr_{\varepsilon} \vc u_{\varepsilon}\cdot \Grade \vc U
\right)
\cdot \left(\vc U-\vc u_{\varepsilon}\right)
-p(\vr_\ep,\vt_{\varepsilon})\ \Dive \vc U
+\tn{S}(\vt_{\varepsilon},\Grade\vu_\ep):\Grade \vc U
\right)\ \dx\ \dt.
\eF

Using  finally as test function in (\ref{m17}) $\varphi=\Theta$, we get
\[
\int_{\Omega}  \vr_{0,\ep} s(\vr_{0,\varepsilon},\vt_{0,\ep}) \Theta (0, \cdot)\ \dx
-\int_{\Omega}  \vr_\ep s(\vr_{\varepsilon},\vt_\ep) \Theta (\tau, \cdot)\ \dx
\]
\[
+\int_0^{\tau}\int_{\Omega}
\frac{\Theta}{\vt_{\varepsilon}} \left( \tn{S}(\vt_{\varepsilon},\Grade \vu_\ep) : \Grade \vu_{\varepsilon}
 - \frac{\vc{q}(\vt_\ep,\Grade \vt_{\varepsilon}) \cdot \Grade \vt_{\varepsilon} }{\vt_{\varepsilon}} \right)\ \dx\ \dt
\]
\[
\leq
-\int_0^{\tau} \int_{\Omega}
\Big(
  \vr_{\varepsilon} s(\vr_{\varepsilon},\vt_\ep) \partial_t \Theta
+\vr_{\varepsilon} s(\vr_{\varepsilon},\vt_\ep)\vu_{\varepsilon}
\cdot \Grade \Theta \Big)\ \dx \ \dt
\]
\bFormula{entro1}
-\int_0^{\tau} \int_{\Omega}
\frac{ {\vc q}(\vt_\ep,\Grade \vt_{\varepsilon})  }{\vt_{\varepsilon}}\cdot \Grade \Theta  \ \dx \ \dt.
\eF

>From (\ref{enertot1}) and (\ref{entro1}) we get
\[
\int_{\Omega}
\left(\frac{1}{2}\ \vr_{\varepsilon} |\vc u_{\varepsilon}-\vc U|^2
+\vr_{\varepsilon}e(\vr_{\varepsilon},\vt_\ep)
-\vr_\ep s(\vr_{\varepsilon},\vt_\ep) \Theta \right)(\tau,\cdot)\ \dx
\]
\[
+\int_0^{\tau}\int_{\Omega}
\frac{\Theta}{\vt_{\varepsilon}} \left( \tn{S}(\vt_{\varepsilon},\Grade\vu_\ep) : \Grade \vu_{\varepsilon}
 - \frac{\vc{q}(\vt_\ep,\Grade \vt_{\varepsilon}) \cdot \Grade \vt_{\varepsilon} }{\vt_{\varepsilon}} \right)\ \dx\ \dt
\]
\[
\leq
\int_{\Omega}
\left(\frac{1}{2} \vr_{0,\varepsilon} |\vc u_{0,\varepsilon}- \vc U(0,\cdot)|^2
+\vr_{0,\varepsilon}e(\vr_{0,\varepsilon},\vt_{0,\ep})
-\vr_{0,\ep} s(\vr_{0,\varepsilon},\vt_{0,\ep}) \Theta (0, \cdot)
\right)\ \dx
\]
\[
+\int_0^{\tau}\int_{\Omega}
\left(
\left(
\vr_{\varepsilon}  \partial_t\vc U
+\vr_{\varepsilon} \vc u_{\varepsilon}\cdot \Grade \vc U
\right)
\cdot \left(\vc U-\vc u_{\varepsilon}\right)
-p(\vr_{\varepsilon},\vt_\ep)\ \Dive \vc U
+\tn{S}(\vr_{\varepsilon},\Grade\vu):\Grade \vc U
\right)\ \dx\ \dt
\]
\bFormula{enerentro}
-\int_0^{\tau} \int_{\Omega}
\Big(  \vr_{\varepsilon} s(\vr_{\varepsilon},\vt_\ep) \partial_t \Theta
+\vr_{\varepsilon} s(\vr_{\varepsilon},\vt_\ep)\vu_{\varepsilon}
\cdot \Grade \Theta +\frac{ {\vc q}_{\varepsilon} }{\vt_{\varepsilon}}\cdot \Grade \Theta
\Big) \ \dx \ \dt
\eF
\[
+\int_0^{\tau} \int_{\Omega}
\left(\ep^{-2\beta}
\vr_{\varepsilon}\vc{E}\cdot (\vc u _{\varepsilon}-\vc U)
 -\vr_{\varepsilon}(\vc\chi \times \vc u_{\varepsilon})\cdot(\vc u_{\varepsilon}-\vc U) + \vr_{\varepsilon}\Grade |\vc \chi \times \vc x|^2\cdot (\vc u_{\varepsilon} - \vc U) \right)\ \dx \ \dt.
\]

Using as test function in (\ref{m12}) with $b \equiv 0$ the function $\phi=\partial_{r}H_{\Theta}(r,\Theta)$, we get
\[
\int_{\Omega} \vr_{\varepsilon}\partial_r H_{\Theta}(r,\Theta)(\tau,\cdot)\ \dx
-\int_{\Omega} \vr_{0,\varepsilon}\partial_r H_{\Theta(0,\cdot)}\left(r(0,\cdot),\Theta(0,\cdot)\right)\ \dx
\]
\bFormula{mass2}
=\int_0^{\tau}\int_{\Omega}
 \Big(
\vr_{\varepsilon}\partial_t\Big(\partial_r H_{\Theta}(r,\Theta)\Big)
+\vr_{\varepsilon}\vc u_{\varepsilon}\cdot \Grade \Big(\partial_r H_{\Theta}(r,\Theta)\Big)
\Big)\ \dx\ \dt.
\eF

Formulas (\ref{enerentro}) and (\ref{mass2}) yield
\[
\int_{\Omega}
\left(
\frac{1}{2}\ \vr_{\varepsilon} |\vc u_{\varepsilon}-\vc U|^2
+H_{\Theta}(\vr_{\varepsilon},\vt_{\varepsilon})
-H_{\Theta}(r,\Theta)
-\partial_r H_{\Theta}(r,\Theta)(\vr_{\varepsilon}-r))
\right)(\tau,\cdot)\ \dx
\]
\[
+\int_0^{\tau}\int_{\Omega}
\frac{\Theta}{\vt_{\varepsilon}} \left( \tn{S}(\vt_\ep,\Grade \vu_{\varepsilon}) : \Grade \vu_{\varepsilon}
 - \frac{\vc{q}_{\varepsilon} \cdot \Grad \vt_{\varepsilon} }{\vt_{\varepsilon}} \right)\ \dx\ \dt
\]
\[
\leq
\int_{\Omega}
\frac{1}{2}\ \vr_{0,\varepsilon} |\vc u_{0,\varepsilon}- \vc U(0,\cdot)|^2\ \dx
\]
\[
+\int_{\Omega}
\big(
H_{\Theta(0,\cdot)}(\vr_{0,\varepsilon},\vt_{0,\varepsilon})
-H_{\Theta(0,\cdot)}(r(0,\cdot),\Theta(0,\cdot))
\]
\[
-\partial_r H_{\Theta(0,\cdot)}(r(0,\cdot),\Theta(0,\cdot))
(\vr_{0,\varepsilon}-r(0,\cdot))
\big)\ \dx
\]
\[
+\int_0^{\tau}\int_{\Omega}
\left(
\left(
\vr_{\varepsilon}  \partial_t\vc U
+\vr_{\varepsilon} \vc u_{\varepsilon}\cdot \Grade \vc U
\right)
\cdot \left(\vc U-\vc u_{\varepsilon}\right)
-p(\vt_{\varepsilon},\vr_\ep)\ \Dive \vc U
+\tn{S}(\vt_{\varepsilon},\Grade \vu_\ep):\Grade \vc U
\right)\ \dx\ \dt
\]
\[
-\int_0^{\tau} \int_{\Omega}
\left(
  \vr_{\varepsilon} s(\vr_{\varepsilon},\vt_\ep) \partial_t \Theta
+\vr_{\varepsilon} s(\vr_{\varepsilon},\vt)\vu_{\varepsilon}
\cdot \Grade \Theta
+\frac{ {\vc q}(\vt_\ep,\Grade \vt_{\varepsilon})  }{\vt_{\varepsilon}}\cdot \Grade \Theta
\right) \ \dx \ \dt
\]
\[
-\int_0^{\tau}\int_{\Omega}
 \Big(
\vr_{\varepsilon}\partial_t\Big(\partial_r H_{\Theta}(r,\Theta)\Big)
+\vr_{\varepsilon}\vc u_{\varepsilon}\cdot \Grade \Big(\partial_r H_{\Theta}(r,\Theta)\Big)
\Big)\ \dx\ \dt
\]
\[
 +
\int_0^{\tau}\int_{\Omega}\left(\ep^{-2\beta} \vr_{\varepsilon}\vc{E} \cdot(\vc u _{\varepsilon}-\vc U)
 -\vr_{\varepsilon}(\vc \chi \times \vc u_{\varepsilon})\cdot (\vc u_{\varepsilon}-\vc U)+ \vr_{\varepsilon}\Grade |\vc \chi \times \vc x|^2 \cdot (\vc u_{\varepsilon} - \vc U)\right)\ \dx \ \dt
\]
\bFormula{enerentro1}
+\int_0^{\tau}\int_{\Omega}
 \partial_t
\Big(r\partial_r H_{\Theta}(r,\Theta)-H_{\Theta}(r,\Theta)\Big)
\ \dx\ \dt.
\eF

Observing finally that for $D=\partial_t$ or $D=\Grad$ one has
\[
D\partial_r H_{\Theta}(r,\Theta)
=-s(r,\Theta)D\Theta
-r\partial_r s(r,\Theta)D\Theta
+\partial_{rr}^2H_{\Theta}(r,\Theta)D r
+\partial_{\Theta,r}^2H_{\Theta}(r,\Theta)D\Theta,
\]
and using the thermodynamic relations (\ref{i5}) (Maxwell's relation) and (\ref{m4a}) (Gibbs' relation),
the following relative entropy inequality holds

\bFrame{
\begin{equation}\label{REI}
\int_{\Omega}
\left(
\frac{1}{2}\ \vr_{\varepsilon} |\vc u_{\varepsilon}-\vc U|^2
+{\mathcal E}\left(\vr_{\varepsilon},\vt_{\varepsilon}|r,\Theta\right)
 \right)(\tau,\cdot)\ \dx
\end{equation}
\[
+\int_0^{\tau}\int_{\Omega}
 \frac{\Theta}{\vt_{\varepsilon}} \left( \tn{S}(\vt_\ep,\Grade \vu_\ep) : \Grade \vu_{\varepsilon}
 - \frac{\vc{q}_{\varepsilon} \cdot \Grad \vt_{\varepsilon} }{\vt_{\varepsilon}} \right)\ \dx\ \dt
\]
\[
\leq
\int_{\Omega}
\left(\frac{1}{2}
 \vr_{0,\varepsilon} |\vc u_{0,\varepsilon}- \vc U(0,\cdot)|^2
+{\mathcal E}
(
\vr_{0,\varepsilon},\vt_{0,\varepsilon}|r(0,\cdot),\Theta(0,\cdot)
)\right)\ \dx
\]
\[
+\int_0^{\tau}\int_{\Omega}
\vr_{\varepsilon} (\vc u_{\varepsilon}-\vc U) \cdot \Grade \vc U
\cdot \left(\vc U-\vc u_{\varepsilon}\right)\ \dx\ \dt
\]
\[
+\int_0^{\tau}\int_{\Omega}
\vr_{\varepsilon}
\left(
s(\vr_{\varepsilon},\vt_\ep)-s(r,\Theta)
\right)
\left(\vc U-\vc u_{\varepsilon}\right)
\cdot\Grade \Theta\ \dx\ \dt
\]
\[
+\int_0^{\tau}\int_{\Omega}
\left(
\vr_{\varepsilon}
\left(
 \partial_t\vc U
+\vc U\cdot \Grade \vc U
\right)
\cdot \left(\vc U-\vc u_{\varepsilon}\right)
\right)\ \dx\ \dt
\]
\[
-\int_0^{\tau}\int_{\Omega}
\left(
p(\vr_{\varepsilon},\vt_\ep)\ \Dive \vc U
-\tn{S}(\vt_\ep,\Grade\vu_{\varepsilon}):\Grade \vc U
\right)\ \dx\ \dt
\]
\[
-\int_0^{\tau}\int_{\Omega}
 \Big(
\vr_{\varepsilon}\left(s(\vr_{\varepsilon},\vt_\ep)-s(r,\Theta)\right)
\partial_t\Theta\Big)\ \dx\ \dt
\]
\[
-\int_0^{\tau} \int_{\Omega} \vr_{\varepsilon}\left(s(\vr_{\varepsilon},\vt_\ep)-s(r,\Theta)\right)
\vc U\cdot  \Grade \Theta \ \dx\ \dt
\]
\[
-\int_0^{\tau}\int_{\Omega}\frac{\vc q(\vt_\ep,\Grade \vt_{\varepsilon}) }{\vt_{\varepsilon}}\cdot  \Grade \Theta
\ \dx\ \dt
\]
\[
+\int_0^{\tau}\int_{\Omega}
\Big(
\left(
1-\frac{ \vr_{\varepsilon}}{r}
\right)
\partial_t p(r,\Theta)
-\frac{ \vr_{\varepsilon}}{r}\ \vc u_{\varepsilon}\cdot \Grade  p(r,\Theta)
\Big)
\ \dx\ \dt
\]
\[ +
\int_0^{\tau}\int_{\Omega}\left(\ep^{-2\beta}\vr_{\varepsilon}\vc{E} \cdot (\vc u _{\varepsilon}-\vc U)
 -\vr_{\varepsilon}(\vc \chi \times \vc u_{\varepsilon})\cdot(\vc u_{\varepsilon}-\vc U) +  \vr_{\varepsilon}\Grade |\vc \chi \times \vc x|^2 \cdot (\vc u_{\varepsilon} - \vc U)\right)\ \dx \ \dt.
 \]
}

\section{ Proof of Theorem \ref{main}}
\label{pr}

\subsection{Preliminaries}
\label{prI}

It is easy to verify that
\bFormula{stress}
\tn{S}(\vt,\nabla _{\epsilon} \vc v):\nabla _{\epsilon} \vc v= \Big(\eta(\vt)-\frac 23 \mu(\vt)\Big)|\Dive \vc v|^2 +
 \mu(\vt) (|\Grade \vc v|^2 + \Grade \vc v : (\Grade \vc v)^T)
\eF
for any $\vc{v} \in W^{1,2}(\Omega;\R^3)$. As for any $\vc{v} \in W^{1,2}_{0,n}(\Omega;\R^3)$
$$
\intO{ \Grade \vc{v} : (\Grade \vc{v})^T} = \intO{(\Dive \vc{v})^2},
$$
we have for
$\vc{v} \in W^{1,2}_{0,n}(\Omega;\R^3)$
\bFormula{Korn1}
\intO{\tn{S}(\vt,\nabla _{\epsilon} \vc v):\nabla _{\epsilon} \vc v} \geq C \|\vc{v}\|_{W^{1,2}(\Omega;\RRR^3)}^2,
\eF
\bFormula{Korn2}
\intO{\frac {1}{\vt}\tn{S}(\vt,\nabla _{\epsilon} \vc v):\nabla _{\epsilon} \vc v} \geq C \|\vc{v}\|_{W^{1,2}(\Omega;\RRR^3)}^2,
\eF
provided $\mu$ fulfills (\ref{m7}), $\eta \equiv 0$, $\epsilon \leq 1$ and $\vt>0$ in $(0,T)\times \Omega$. Indeed, easily
\bFormula{Korn1a}
\intO{\tn{S}_h(\Theta,\Gradh \vc V):\Gradh \vc V} \geq C \|\vc{V}\|_{W^{1,2}(\omega;\RRR^2)}^2,
\eF
\bFormula{Korn2a}
\intO{\frac {1}{\Theta}\tn{S}_h(\Theta,\nabla _{h} \vc V):\nabla _{h} \vc V} \geq C \|\vc{V}\|_{W^{1,2}(\omega;\RRR^2)}^2
\eF
for any $\vc{V} \in W^{1,2}_0(\omega;\R^2)$ and $\Theta > 0$ in $(0,T)\times \omega$.

Moreover, note that we also have the Poincar\'e inequality in the form
\bFormula{stress2}
\|\vc V \|_{L^2(\omega;\RRR^2)} \leq c\|\nabla _h \vc V\|_{L^2(\omega;\RRR^{2\times 2})}
\eF
for any $\vc{V} \in W^{1,2}_0(\omega;\R^2)$.

Due to the energy equality (\ref{m18}) combined with the entropy inequality (\ref{m17}) and Korn's inequality in the form (\ref{Korn2}), we have that the sequence $(\vr_\ep,\vu_\ep,\vt_\ep)$ is bounded in the following spaces
\bFormula{est1}
\|\vr_\ep\|_{L^\infty(0,T;L^\gamma(\Omega))} + \|\sqrt{\vr_\ep} \vu_\ep\|_{L^\infty(0,T;L^2(\Omega;\RRR^3))} + \|\vu\|_{L^2(0,T;W^{1,2}(\Omega;\RRR^3))} + \|\Grade \vt_\ep\|_{L^2(0,T;L^2(\Omega;\RRR^3))}
\eF
\[ +\|\vt_\ep\|_{L^\infty(0,T;L^4(\Omega))} + \|\vt_\ep\|_{L^3(0,T;L^9(\Omega))} \leq C
\]
with the constant $C$ independent of $\epsilon$. These estimates hold if $\gamma \geq \frac{12}{5}$ (if $\alpha=1$) or under the assumptions on $g$ from Theorem \ref{main} (if $\alpha =0$), for any $\gamma > \frac 32$. Note that the limit on $\gamma$ comes from the gravitational potential, as
\[
\Big\|\Grade \int_\Omega \frac{\vr(y)}{(\sqrt{(x_h-y_h)^2 + \ep^2(x_3-y_3)})^3}\ \dy \Big\|_{L^p(\Omega;\RRR^3)} \leq C \|\vr\|_{L^p(\Omega)}
\]
for $1<p<\infty$, with $C$ independent of $\ep$.

We recall the necessary definitions of essential and residual sets from \cite{FEINOV}.
We take four positive numbers $0<\underline{\vr}\leq \overline{\vr}<\infty$, $0<\underline{\theta}\leq \overline{\theta}<\infty$.
The essential and residual subsets of $\Omega$ are defined for a.e. $t\in (0,T)$ as follows:
\bFormula{essres}
{\mathcal O}_{ess}^{\ep}(t)=\Big\{ x\in\Omega\ | \frac{1}{2}\ \underline{\vr}\leq \vr_{\ep}(t,x)\leq 2\overline{\vr}, \, \frac{1}{2}\ \underline{\vt}\leq \vt_{\ep}(t,x)\leq 2\overline{\vt} \Big\},\ \ \
{\mathcal O}_{res}^{\ep}(t)=\Omega\backslash {\mathcal O}_{ess}^{\ep}(t).
\eF
For any function $h$ defined for a.e. $(t,x)\in (0,T)\times\Omega$, we write
\bFormula{essres2}
[h]_{ess}(t,x)=h(t,x) {\mathbf 1}_{{\mathcal O}_{ess}^{\ep}(t)}(x),\ \ \ [h]_{res}(t,x)=h(t,x) {\mathbf 1}_{{\mathcal O}_{res}^{\ep}(t)}(x).
\eF


Now, using \cite{FEINOV}, we have the following properties of the  Helmholtz function

\bLemma{l1}
Let $\overline{\vr}>0$ and $\overline{\vt}>0$ be two given constants and let $H_{\overline{\vt}}(\vr,\vt)=\vr e-\overline{\vt}\vr s$. Let ${\mathcal O}_{ess}(t)$ be defined as above with $\overline{\vr}=\underline{\vr}$ and $\overline{\vt}=\underline{\vt}$.
Then there exist positive constants $C_j=C_j(\overline{\vr},\overline{\vt})$ for $j=1,\dots,4$ such that
\[
C_1\left(|\vr-\overline{\vr}|^2+|\vt-\overline{\vt}|^2\right)
\leq
H_{\overline{\vt}}(\vr,\vt)
-(\vr-\overline{\vr})\partial_{\vr}H_{\overline{\vt}}(\overline{\vr},\overline{\vt})
-H_{\overline{\vt}}(\overline{\vr},\overline{\vt})
\]
\bFormula{ess1}
\leq
C_2\left(|\vr-\overline{\vr}|^2+|\vt-\overline{\vt}|^2\right),
\eF
for all $(\vr,\vt)\in{\mathcal O}_{ess}(t)$,
\[
H_{\overline{\vt}}(\vr,\vt)
-(\vr-\overline{\vr})\partial_{\vr}H_{\overline{\vt}}(\overline{\vr},\overline{\vt})
-H_{\overline{\vt}}(\overline{\vr},\overline{\vt})
\]
\bFormula{res1}
\geq
\inf_{(\tilde\vr,\tilde\vt) \in{\partial \mathcal O}_{ess}}
\left\{H_{\overline{\vt}}(\tilde\vr,\tilde\vt)
-(\tilde\vr-\overline{\vr})\partial_{\vr}H_{\overline{\vt}}(\overline{\vr},\overline{\vt})
-H_{\overline{\vt}}(\overline{\vr},\overline{\vt})\right\}
= C_3,
\eF
for all $(\vr,\vt)\in{\mathcal O}_{res}(t)$,
\bFormula{res2}
H_{\overline{\vt}}(\vr,\vt)
-(\vr-\overline{\vr})\partial_{\vr}H_{\overline{\vt}}(\overline{\vr},\overline{\vt})
-H_{\overline{\vt}}(\overline{\vr},\overline{\vt})
\geq
C_4\left( \vr e(\vr,\vt)+\vr |s(\vr,\vt)|\right),
\eF
for all $(\vr,\vt)\in{\mathcal O}_{res}(t)$,
\eL
\begin{proof}
See  \cite[Lemma 5.1]{FEINOV}.
\end{proof}

As a consequence we get the following lemma
\bLemma{res}
\label{res}
There exists a constant $C=C(\underline{\vr},\overline{\vr},\underline{\vt}, \overline{\vt})>0$ such that for all $\vr\in[0,\infty)$, $r\in[\underline{\vr}/2,2\overline{\vr}]$, $\vt \in (0,\infty)$ and $\Theta \in [\underline{\vt}/2,2\overline{\vt}]$
\[
{\mathcal E}\left(\vr,\vt|r,\Theta \right)
 (\tau,\cdot) \geq
\]
\[
 C(\underline{\vr},\overline{\vr},\underline{\vt}, \overline{\vt})\left({\mathbf 1}_{{\mathcal O}_{res}}+\vr^{\gamma}{\mathbf 1}_{{\mathcal O}_{res}}+\vt^{4}{\mathbf 1}_{{\mathcal O}_{res}}+
(\vr-r)^2{\mathbf 1}_{{\mathcal O}_{ess}}+ (\vt-\Theta)^2{\mathbf 1}_{{\mathcal O}_{ess}}\right).
\]
\eL

The lemma yields the lower bound of the relative entropy functional
\bFormula{lb}
{\mathcal I}(\vr_{\ep}, \vu_{\ep},\vt_{\ep}| r, \vec U, \Theta ) \geq \
\eF
\[
C(\underline{\vr},\overline{\vr},\underline{\vt},\overline{\vt}))
\int_{\Omega}\left(\vr_\ep |\vu_\ep-\vc{U}|^2 + {\mathbf 1}_{res}+[\vr_{\ep}^{\gamma}]_{res}+[\vr_\ep-r]^2_{ess}+ [\vt_{\ep}^{4}]_{res}+[\vt_\ep-\Theta]^2_{ess}\right)\ \dx.
\]




\subsection{Estimates of the remainder}
\label{prII}

We now return back to (\ref{REI}), where we assume that $(r,\vc{U},\Theta)$, $\vc{U}=(\vc{V},0)$, is such that $(r,\vc{V},\Theta)$ solves the 2-D Navier--Stokes--Fourier--Poisson system (\ref{t1}--\ref{t5}) in $(0,T)\times \omega$, see Proposition \ref{pro1}. We further use the following. In order to integrate over $\Omega$, we assume that the functions defined only on $\omega$ are extended being constant in $x_3$ for $0\leq x_3\leq 1$. Moreover, we write $\vc{U}$ instead of $\vc{V}$ in the situations, when we need to use a vector field with three components. Similarly, we write in such situations $\tn{S}(\Theta,\Grade U)$ instead of $\tn{S}_h (\Theta, \Gradh \vc{V})$. Recall that these functions are sufficiently regular and fulfill the corresponding boundary conditions, therefore they can be used as test functions in (\ref{m12}--\ref{m17}).

We denote $\underline{\vr} = \inf_{(0,T)\times \Omega} r$, $\overline{\vr} = \sup_{(0,T)\times \Omega} r$, $\underline{\vt} = \inf_{(0,T)\times \Omega} \Theta$ and $\overline{\vt} = \sup_{(0,T)\times \Omega} \Theta$ and use these numbers in order to define the essential and residual sets in (\ref{essres}). As above,  $(\vr_{\ep},\vu_{\ep},\vt_{\ep} )$ denotes the solution to the primitive system.   We now rearrange the remainder, i.e. the terms on the right-hand side of (\ref{REI}) which are integrated over time and space, as follows
\begin{equation}
\label{Re1}
\begin{array}{c}
\displaystyle \int_0^{\tau}{\mathcal R}(\vr_{\ep}, \vu_{\ep},\vt_{\ep},  r, \vec V,\Theta )\ \dt
 \\ \\
\displaystyle =\int_0^{\tau}\int_{\Omega}
\vr_{\ep}(\vu_{\ep}-\vc U)\cdot\Grade \vc U\cdot(\vc U-\vu_{\ep})\ \dx \ \dt \\[10pt]
\displaystyle  + \int_0^{\tau}\int_{\Omega}
\vr_{\varepsilon} \left(
s(\vr_{\varepsilon},\vt_\ep)-s(r,\Theta)
\right)
\left(\vc U-\vc u_{\varepsilon}\right)
\cdot\Grade \Theta\ \dx\ \dt
\\[10pt]
\displaystyle +\int_0^{\tau}\int_{\Omega}
\vr_{\ep}\left(\partial_t\vc U+\vc U\cdot\Grade \vc U\right)\cdot(\vc U- \vu_{\ep})
\ \dx\ \dt  \\[10pt]
\displaystyle -\int_0^{\tau}\int_{\Omega}
 \Big(
\vr_{\varepsilon}\left(s(\vr_{\varepsilon},\vt_\ep)-s(r,\Theta)\right)
\partial_t\Theta\Big)\ \dx\ \dt
\\[10pt]
\displaystyle -\int_0^{\tau} \int_{\Omega} \vr_{\varepsilon}\left(s(\vr_{\varepsilon},\vt_\ep)-s(r,\Theta)\right)
\vc V\cdot  \Gradh \Theta \ \dx\ \dt \\[10pt]
\displaystyle +\int_0^{\tau}\int_{\Omega}
\Big(
\left(
1-\frac{ \vr_{\varepsilon}}{r}
\right)
\partial_t p(r,\Theta)
-\frac{ \vr_{\varepsilon}}{r}\ \vc u_{\varepsilon}\cdot \Grade  p(r,\Theta)
\Big)
\ \dx\ \dt  \\[10pt]
\displaystyle +\int_0^{\tau}\int_{\Omega}
\vr_{\ep}(\vec \chi\times\vu_{\ep})\cdot(\vc U-\vu_{\ep})
\ \dx\ \dt
\displaystyle - \int_0^{\tau}\int_{\Omega}
\vr_{\ep}\Grade |\vec \chi\times\vec x|^2\cdot(\vc U-\vu_{\ep})
\ \dx\ \dt \\[10pt]
\displaystyle  - \int_0^{\tau}\int_{\Omega} \ep^{-2\beta}
\vr_{\ep}  \vc{E} \cdot(\vc U-\vu_{\ep})
\ \dx\ \dt
\displaystyle -\int_0^{\tau}\int_{\Omega}\frac{\vc q(\vt_\ep,\Grade\vt_\ep)}{\vt_{\varepsilon}}\cdot  \Grade \Theta
\ \dx\ \dt \\[10pt]
\displaystyle
-\int_0^{\tau}\int_{\Omega}
\left(
p(\vr_{\varepsilon},\vt_\ep)\ \Divh \vc V
-\tn{S}(\vt_{\varepsilon},\Grade \vu_\ep):\Grade \vc U
\right)\ \dx\ \dt
=: \int_0^\tau \Big(\sum_{j=1}^{11} R_j\Big) \ \dt.
\end{array}
\end{equation}

Let us now estimate all  terms in ${\mathcal R}(\vr_\ep, \vu_\ep, \vt_\ep,  r, \vec V,\Theta)$.

We have
\[ R_1=
 \int_{\Omega}
\vr_{\ep}(\vu_{\ep}-\vc U)\cdot\Grade \vc U\cdot(\vc U-\vu_{\ep})\ \dx
\leq
C D(t){\mathcal I}(\vr_{\ep}, \vu_{\ep},\vt_{\ep}| r, \vec U, \Theta)
\]
with
\[
D=\|\Gradh \vc V\|_{L^{\infty}(\Omega;\RRR^{2\times 2})}\in L^1(0,T),
\]
due to the assumption on the smoothness of $\vc{V}$. Next


\begin{equation}
\begin{array}{rcl}
|R_2|
&\leq &
\displaystyle \left|
\int_{\Omega}
\vr_{\varepsilon}
\left(
s(\vr_{\varepsilon},\vt_\ep)-s(r,\Theta)
\right)
(\vc U-\vc u_{\varepsilon})
\cdot\Grade \Theta \ \dx
\right|
\end{array}
\end{equation}

\[
\leq \| \Gradh \Theta \|_{L^{\infty}(\Omega;\RRR^2)}
\Big[
2\overline{\rho}
\int_{\Omega}
\big|
\big[
s(\vr_{\varepsilon},\vt_{\varepsilon})-s(r,\Theta)
\big]_{ess}
\big| \,
|\vc U-\vc u_{\varepsilon}|
\ \dx
\]
\[
+  \int_{\Omega}
\big|
\big[
\vr_{\varepsilon}
(
s(\vr_{\varepsilon},\vt_{\varepsilon})-s(r,\Theta)
)
\big]_{res}
\big|\,
|\vc U-\vc u_{\varepsilon}|
\ \dx
\Big].
\]

Lemma \ref{res} together with the properties of the entropy (\ref{Ent1}--\ref{Ent2}) yields
\[
\int_{\Omega}
\big|
\big[
s(\vr_{\varepsilon},\vt_{\varepsilon})-s(r,\Theta)
\big]_{ess}
\big|\,
|\vc U-\vc u_{\varepsilon}|
\ \dx
\leq
\delta\|\vc U-\vc u_{\varepsilon}\|^2_{L^2(\Omega;\RRR^3)}
+C(\delta)
\int_{\Omega}
{\mathcal E}(\vr_{\varepsilon},\vt_{\varepsilon}|r,\Theta)\ \dx,
\]
for $\delta>0$, and
\[
\int_{\Omega}
\big|
\big[
\vr_{\varepsilon}
(
s(\vr_{\varepsilon},\vt_{\varepsilon})-s(r,\Theta)
)
\big]_{res}
\big|
\, |\vc U-\vc u_{\varepsilon}|
\ \dx
 \]
\[
\leq \delta\|\vc U-\vc u_{\varepsilon}\|^2_{L^6(\Omega;\RRR^3)}
+C(\delta)
\|
\big[
\vr_{\varepsilon}
(
s(\vr_{\varepsilon},\vt_{\varepsilon})-s(r,\Theta)
)
\big]_{res}
\|^2_{L^{6/5}(\Omega)}.
\]
Using again properties of the entropy (\ref{Ent1}--\ref{Ent2}) together with the fact that the mapping  $t\mapsto \int_{\Omega}
{\mathcal E}\left(\vr_{\varepsilon},\vt_{\varepsilon}|r,\Theta\right)\ \dx\in L^{\infty}(0,T)$, we conclude that
\[
\left\|
\left[
\vr_{\varepsilon}
\left(
s(\vr_{\varepsilon},\vt_{\varepsilon})-s(r,\Theta)
\right)
\right]_{res}
\right\|^2_{L^{6/5}(\Omega)}
\leq C
\int_{\Omega}{\mathcal E}\left(\vr_{\varepsilon},\vt_{\varepsilon}|r,\theta\right)\ \dx.
\]
So finally, we end up with
\[
|R_2|
\leq
\delta\|\vc U-\vc u_{\varepsilon}\|^2_{W_0^{1,2}(\Omega;\RRR^3)}
+C(\delta;r,\vc V,\Theta)
\int_{\Omega}{\mathcal E}\left(\vr_{\varepsilon},\vt_{\varepsilon}|r,\Theta\right)\ \dx.
\]
Next, using the fact that $(r,\vc{V},\Theta)$ solve the 2-D Navier--Stokes--Fourier--Poisson system,
\[
R_3=\int_{\Omega}
\vr_{\varepsilon}
(
 \partial_t\vc U
+\vc U\cdot \Gradh \vc U
)
\cdot (\vc U-\vc u_{\varepsilon})
\ \dx
= R_{3,1}+R_{3,2},
\]
where
\[ R_{3,1}= \int_{\Omega}
\frac{\vr_{\varepsilon}}{r} (\vc U-\vc u_{\varepsilon})\cdot
(
\Dive \tn{S}(\Theta,\Grade \vc U)-\Grade p(r,\Theta)
)\ \dx,
\]
\[ R_{3,2}= \int_{\Omega} \vr_{\varepsilon}(\vc U-\vc u_{\varepsilon})\cdot (- (\vc \chi \times \vc U)
+  \Grade |\vc \chi\times \vc x|^2 +  \Gradh \phi_h)\ \dx =\sum _{i=1}^3 K_i.
\]
We write
\[R_{3,1}
=\int_{\Omega}
\frac{\vr_{\varepsilon}-r}{r}(\vc U-\vc u_{\varepsilon})\cdot
(
\Dive \tn{S}(\Theta,\Grade \vc U)
-\Grade p(r,\Theta)
)\ \dx\]
\[
+\int_{\Omega}
 (\vc U-\vc u_{\varepsilon})\cdot
(
\Dive \tn{S}(\Theta,\Grade \vc U)
-\Grade p(r,\Theta)
)\ \dx.
\]
We now divide the first term in $R_{3,1}$ into the integral over the essential and the residual set and estimate them similarly as $R_2$. It reads
\[
\Big|
\int_{\Omega}
\frac{\vr_{\varepsilon}-r}{r}(\vc U-\vc u_{\varepsilon})\cdot
(
\Dive \tn{S}(\Theta,\Grade \vc U)
-\Grade p(r,\Theta)
)\ \dx\Big|
\]\[\leq
C(\delta;r,\vc V,\Theta)
\|\left[\vr_{\varepsilon}-r\right]_{ess}\|_{L^2(\Omega)}^2
+
\delta\|\vc U-\vc u_{\varepsilon}\|^2_{L^2(\Omega;\RRR^3)}
\]
\[
+
C(\delta;r,\vc V,\Theta)
\left(
\|\left[\vr_{\varepsilon}\right]_{res}\|_{L^{6/5}(\Omega)}^2
+\|\left[1\right]_{res}\|_{L^{6/5}(\Omega)}^2
\right)
+
\delta\|\vc U-\vc u_{\varepsilon}\|^2_{L^{6}(\Omega;\RRR^3)}.
\]
Integrating by parts in the second integral of $R_{3,1}$, we have
\[
\int_{\Omega}
 (\vc U-\vc u_{\varepsilon})\cdot
(
\Dive \tn{S}(\Theta,\Grade \vc U)
-\Grade p(r,\Theta)
\Big)\ \dx
\] \[
=-\int_{\Omega}
(
 \tn{S}(\Theta,\Grade \vc U):\Grade (\vc U-\vc u_{\varepsilon})
- p(r,\Theta)\ \Dive (\vc U-\vc u_{\varepsilon}))\ \dx.
\]
Hence we conclude
\[
R_{3,1}
\leq
\int_{\Omega}
\Big(
  p(r,\Theta)\ \Divh (\vc V-\vc u_{\varepsilon})-\tn{S}(\Theta,\Grade \vc U):\Grade (\vc U-\vc u_{\varepsilon})
\Big)\ \dx
+\delta\|\vc U-\vc u_{\varepsilon}\|^2_{W^{1,2}(\Omega;\RRR^3)}
\]
\[
+C(\delta;r,\vc V,\Theta)
\int_{\Omega}{\mathcal E}\left(\vr_{\varepsilon},\vt_{\varepsilon}|r,\Theta\right)\ \dx,
\]
for any $\delta>0$.

The terms $K_1$--$K_3$ will be treated below, in combination with $R_7$--$R_9$.
%

\[
R_4=
 -\int_{\Omega}
  \vr_{\varepsilon} \big( s(\vr_{\varepsilon} ,\vt_{\varepsilon})- s(r ,\Theta)\big)\partial_t \Theta
 \ \dx
\] \[
=-\int_{\Omega}
  (\vr_{\varepsilon}-r) \big( s(\vr_{\varepsilon} ,\vt_{\varepsilon})- s(r ,\Theta)\big)\partial_t \Theta
 \ \dx - \int_{\Omega}
  r \big( s(\vr_{\varepsilon} ,\vt_{\varepsilon})- s(r ,\Theta)\big)\partial_t \Theta
 \ \dx.
\]
We rewrite the first term as above, using the essential and residual parts, and estimate them
\[
-\int_{\Omega}
  (\vr_{\varepsilon}-r) \big( s(\vr_{\varepsilon} ,\vt_{\varepsilon})- s(r ,\Theta)\big)\partial_t \Theta
 \ \dx
\]
\[
= -\int_{\Omega}
  (\vr_{\varepsilon}-r) \big[ s(\vr_{\varepsilon} ,\vt_{\varepsilon})- s(r ,\Theta)\big]_{ess}\partial_t \Theta
 \ \dx -\int_{\Omega}
  (\vr_{\varepsilon}-r) \big[ s(\vr_{\varepsilon} ,\vt_{\varepsilon})- s(r ,\Theta)\big]_{res}\partial_t \Theta
 \ \dx
\]
\[ \leq C(\delta;r,\vc V,\Theta)
\int_{\Omega}{\mathcal E}\left(\vr_{\varepsilon},\vt_{\varepsilon}|r,\Theta\right)\ \dx.
\]
Finally
\[
-\int_{\Omega}
  r \big( s(\vr_{\varepsilon} ,\vt_{\varepsilon})- s(r ,\Theta)\big)\partial_t \Theta
 \ \dx
\]
\[
=-\int_{\Omega}
  r \big( s(\vr_{\varepsilon} ,\vt_{\varepsilon}) - s(r ,\Theta)
-\partial_{\vr}s(r ,\Theta)(\vr_{\varepsilon} -r)
-\partial_{\vt}s(r ,\Theta)(\vt_{\varepsilon} -\Theta)
\big)\partial_t \Theta
 \ \dx
\]
\[
-\int_{\Omega}
  r \big(
\partial_{\vr}s(r ,\Theta)(\vr_{\varepsilon} -r)
+\partial_{\vt}s(r ,\Theta)(\vt_{\varepsilon} -\Theta)
\big)\partial_t \Theta
 \ \dx.
\]
The first integral in the right-hand side can be estimated in the same way as before, using the essential and the residual parts, and we end with
\[
R_4\leq
C(\delta;r,\vc V,\Theta)
\int_{\Omega}{\mathcal E}\left(\vr_{\varepsilon},\vt_{\varepsilon}|\vr,\Theta\right)\ \dx
-\int_{\Omega}
  r \big(
\partial_{\vr}s(r ,\Theta)(\vr_{\varepsilon} -r)
+\partial_{\vt}s(r ,\Theta)(\vt_{\varepsilon} -\Theta)
\big)\partial_t \Theta
 \ \dx.
\]
In the estimate of the next term we use exactly the same procedure as for $R_4$ and end up with
\[R_5=
-\int_{\Omega}
 \vr_{\varepsilon}\big( s(\vr_{\varepsilon} ,\vt_{\varepsilon})- s(r ,\Theta)\big)\vc V\cdot  \Gradh \Theta
\ \dx
\leq
C(\delta;r,\vc V,\Theta)
\int_{\Omega}{\mathcal E}\left(\vr_{\varepsilon},\vt_{\varepsilon}|r ,\Theta\right)\ \dx
\]
\[
-\int_{\Omega}
  r \big(
\partial_{\vr}s(r ,\Theta)(\vr_{\varepsilon} -r)
+\partial_{\vt}s(r ,\Theta)(\vt_{\varepsilon} -\Theta)
\big)\vc V\cdot  \Gradh \Theta
 \ \dx.
\]
Furthermore
\[
R_6=\int_{\Omega}
\Big(
\left(
1-\frac{ \vr_{\varepsilon}}{r}
\right)
\partial_t p(r ,\Theta)
-\frac{ \vr_{\varepsilon}}{r}\ \vc u_{\varepsilon}\cdot \Grade  p(r ,\Theta)
\Big)
\ \dx
\]
\[
=
\int_{\Omega}
\Big(
\left(
1-\frac{ \vr_{\varepsilon}}{r}
\right)
\Big(
\partial_t p(r ,\Theta)
+ \vc V\cdot \Gradh  p(r ,\Theta)
\Big)
\ \dx
+
\int_{\Omega}
p(r ,\Theta)\Dive \vc u_{\varepsilon}
\ \dx
\]
\[
+
\int_{\Omega}
\Big(
\left(
1-\frac{ \vr_{\varepsilon}}{r}
\right)
\Grade p(r ,\Theta) \cdot
( \vc u_{\varepsilon}-\vc U)
\ \dx,
\]
and using the same argument as used for $R_2$, we get
\[
\left|
\int_{\Omega}
\Big(
\left(
1-\frac{ \vr_{\varepsilon}}{r}
\right)
\Grade p(r ,\Theta)\cdot
( \vc u_{\varepsilon}-\vc U)
\ \dx
\right|
\]
\[
\leq
C(\delta;r,\vc V,\Theta)
\Big(
\delta\|\vc U-\vc u_{\varepsilon}\|^2_{W^{1,2}(\Omega;\RRR^3)}
+
\int_{\Omega}{\mathcal E}\left(\vr_{\varepsilon},\vt_{\varepsilon}|r ,\Theta\right)\ \dx\Big)
\]
for any $\delta>0$. Hence we end with
\[
R_6
\leq
\int_{\Omega}
\Big(
\Big(
1-\frac{ \vr_{\varepsilon}}{r}
\Big)
\Big(
\partial_t p(r ,\Theta)
+ \vc V\cdot \Gradh  p(r ,\Theta)
\Big)
\ \dx
+
\int_{\Omega}
p(r ,\Theta)\Dive \vc u_{\varepsilon}
\ \dx
\]
\[
+\delta\|\vc U-\vc u_{\varepsilon}\|^2_{W^{1,2}(\Omega;\RRR^3)}
+C(\delta;r,\vc V,\Theta)
\int_{\Omega}{\mathcal E}\left(\vr_{\varepsilon},\vt_{\varepsilon}|r ,\Theta\right)\ \dx.
\]
We can rewrite $R_7$ and $K_1$ as follows (note that $(\vc \chi \times (\vc{U}-\vc{u}_\ep)\cdot (\vc{U}-\vc{u}_\ep) = \vc 0$)
\[
R_7 + K_1=
 \int_{\Omega}
\vr_{\ep}(\vec \chi\times\vc U)\cdot(\vc U-\vu_{\ep})\ \dx -  \int_\Omega \vr_{\ep}(\vec \chi\times\vc U)\cdot(\vc U-\vu_{\ep})
\ \dx = 0.
\]

Similarly, we have
\[
R_8 + K_2 = 0.
\]

Next we consider the gravitational potential. The case $\alpha =0$ (i.e. only the gravitational effect of other objects than the fluid itself is taken into account) is relatively easy. Recall that we assume
\[
\int_{\RR^3} \frac{g(y)y_3}{(\sqrt{|x_h-y_h|^2 + y_3^2})^3} =0.
\]
Therefore we have to show that
\bFormula{G1}
\lim_{\ep \to 0^+} \int_{\Omega} \vr_\ep (\vc{U}-\vc{u}_\ep)\cdot \Big(\int_{\RR^3} g(y)
\eF
\[
\Big[\frac{(x_h-y_h,-y_3)}{(\sqrt{|x_h-y_h|^2 + y_3^2})^3} - \frac{(x_h-y_h,\ep x_3-y_3)}{(\sqrt{|x_h-y_h|^2 + (\ep x_3-y_3)^2})^3}   \Big]\ \dy \Big)\ \dx = 0.
\]
First of all, due to the estimates presented above, it is not difficult to see that (note that to get estimates independent of $\ep$ of the integral over $\R^3$ is easy) it is enough to verify
\[
\lim_{\ep \to 0^+} \int_{\RR^3} \Ov{g}(y) \Big[\frac{(x_h-y_h,-y_3)}{(\sqrt{|x_h-y_h|^2 + y_3^2})^3} - \frac{(x_h-y_h,\ep x_3-y_3)}{(\sqrt{|x_h-y_h|^2 + (\ep x_3-y_3)^2})^3}   \Big]\ \dy  = \vc{0}
\]
for all $x_h \in \omega$, $x_3 \in (0,1)$ and $\Ov{g} \in C^\infty_c(\R^3)$. As
\[\lim_{\ep \to 0^+} \Big(\frac{(x_h-y_h,-y_3)}{(\sqrt{|x_h-y_h|^2 + y_3^2})^3} - \frac{(x_h-y_h,\ep x_3-y_3)}{(\sqrt{|x_h-y_h|^2 + (\ep x_3-y_3)^2})^3}\Big)   = \vc{0}
\]
for a.a. $(x_h,x_3) \in \Omega$, $(y_h,y_3) \in \R^3$, and
\[
\Big| \frac{(x_h-y_h,\ep x_3-y_3)}{(\sqrt{|x_h-y_h|^2 + (\ep x_3-y_3)^2})^3}\Big| \leq \frac{1}{\sqrt{|x_h-y_h|^2 + (\ep x_3-y_3)^2}} \Big| \in L^1_{{\rm loc}}(\R^3) \quad \forall \ep \in [0,1],
\]
the Lebesgue dominated convergence theorem yields the required identity (\ref{G1}).

The case of the self-gravitation ($\alpha =1$) is more complex. Here, we have to show that
\bFormula{G2}
\int_{\Omega} \vr_\ep (\vc{U}-\vc{u}_\ep)\cdot \Big[\int_{\Omega} \frac{\vr_\ep(t,y)(x_h-y_h,\ep (x_3-y_3))}{(\sqrt{|x_h-y_h|^2 + \ep^2(x_3-y_3)^2})^3}\ \dy + \Grade \int_{\omega}  \frac{r(t,y_h)}{|x_h-y_h|} \ \dy_h  \Big]\ \dx
\eF
\[
\leq \delta \| \vc U-\vu_\ep\|_{L^6(\Omega;\RRR^3)} + C(\delta; r, \vc{V},\Theta) \int_\Omega {\mathcal E}(\vr_\ep,\vt_\ep|r,\Theta)\, \dx + H_\ep,
\]
where $H_\ep = o(\ep)$ as $\ep \to 0^+$. The derivative of the integral over $\omega$ with respect to $x_3$ is indeed zero. First of all, for $\gamma \geq \frac {12}{5}$, as in (\ref{est1}), using the decomposition to the essential and the residual set  and proceeding as in the estimates of the remainder above, we can show that it is enough to verify that
\[
\lim_{\ep \to 0^+} \int_{\Omega} r \vc{U}\cdot \Big[\int_{\Omega} \frac{r(t,y_h)(x_h-y_h,\ep^2 (x_3-y_3))}{(\sqrt{|x_h-y_h|^2 + \ep^2(x_3-y_3)^2})^3}\ \dy + \Grade \int_{\omega}  \frac{r(t,y_h)}{|x_h-y_h|} \ \dy_h  \Big]\ \dx = 0.
\]

Again, it is not difficult to verify that (e.g. using the change of the variables to integrate over $\Omega_\ep$) it is enough to show
\[
\lim_{\ep \to 0^+} \Big[\int_{\Omega} \frac{r(t,y_h)\big(x_h-y_h,\ep (x_3-y_3)\big)}{(\sqrt{|x_h-y_h|^2 + \ep^2(x_3-y_3)^2})^3}\ \dy + \Grade \int_{\omega}  \frac{r(t,y_h)}{|x_h-y_h|} \ \dy_h  \Big] = \vc{0}.
\]
First, note that
\[
\Grade \int_{\omega}  \frac{r(t,y_h)}{|x_h-y_h|} \ \dy_h = -{\rm v.p.} \int_\omega \frac{r(t,y_h)(x_h-y_h)}{|x_h-y_h|^{\frac 32}}\ \dy_h,
\]
where v.p. denotes the integral in the principal value sense. Therefore, we have to verify that
\bFormula{G3}
\lim_{\ep \to 0^+} \int_{\Omega} \frac{\ep r(t,y_h) (x_3-y_3)}{(\sqrt{|x_h-y_h|^2 + \ep^2(x_3-y_3)^2})^3}\ \dy = 0,
\eF
and
\bFormula{G4}
\lim_{\ep \to 0^+} \int_{\Omega} \frac{\ep r(t,y_h) (x_h-y_h)}{(\sqrt{|x_h-y_h|^2 + \ep^2(x_3-y_3)^2})^3}\ \dy = {\rm v.p.} \int_\omega \frac{r(t,y_h)(x_h-y_h)}{|x_h-y_h|^3}\ \dy_h.
\eF
Let us fix $x_0 \in \omega$, $\Delta >0$, sufficiently small, and denote $B_\Delta(x_0) = \{x\in \omega;|x-x_0| <\Delta\}$ and $C_\Delta (x_0) = \{x \in \Omega; |x_h-x_0| <\Delta, 0< x_3 <1\}$.

We first consider (\ref{G3}). Let us fix $\delta >0$. Using the change of variables (from $\Omega$ back to $\Omega_\ep$) it is not difficult to see that there exists $\Delta >0$ such that for any $0<\ep \leq 1$, $0<x_3<1$ we have
\[
\Big|\int_{C_\Delta(x_0)} \frac{\ep r(t,y_h) (x_3-y_3)}{(\sqrt{|x_0-y_h|^2 + \ep^2(x_3-y_3)^2})^3}\ \dy\Big| < \delta
\]
and for this $\Delta>0$ there exists $\ep_0>0$ such that for any $0<\ep\leq \ep_0$
\[
\Big|\int_{\Omega \setminus C_\Delta(x_0)} \frac{\ep r(t,y_h) (x_3-y_3)}{(\sqrt{|x_0-y_h|^2 + \ep^2(x_3-y_3)^2})^3}\ \dy\Big| < \delta.
\]
Whence  (\ref{G3}).

In order to prove (\ref{G4}), we proceed similarly. Since $\frac{x_h-y_h}{|x_h-y_h|^3}$ is a singular integral kernel in the sense of the Calder\'on--Zygmund theory, as above, for a fixed $x_0 \in \omega$, $0<x_3<1$ and $\delta>0$ that there there exists $\Delta >0$ such that
\[
\Big|\int_{C_\Delta(x_0)} \frac{r(t,y_h) (x_0-y_h)}{(\sqrt{|x_0-y_h|^2 + \ep^2(x_3-y_3)^2})^3}\ \dy\Big| < \delta,
\]
and
\[
{\rm v.p.} \int_{B_\Delta(x_0)} \frac{r(t,y_h)(x_h-y_h)}{|x_h-y_h|^3}\ \dy_h\Big| <\delta.
\]
For this $\Delta >0$, using that
\[
\frac{1}{(\sqrt{|x_0-y_h|^2 + \ep^2(x_3-y_3)^2})^3} - \frac{1}{|x_0-y_h|^3} \to 0 \quad \text{ as } \ep \to 0^+
\]
for any $y_h \in \omega$, $0<x_3,y_3<1$, except $x_0 = y_h$, we see that for above given $\Delta >0$, there exists $\ep_0>0$ such that for any $0<\ep\leq \ep_0$
\[
\Big|\int_{\Omega\setminus C_\Delta(x_0)} \frac{\ep r(t,y_h) (x_h-y_h)}{(\sqrt{|x_h-y_h|^2 + \ep^2(x_3-y_3)^2})^3}\ \dy - {\rm v.p.} \int_{\omega\setminus B_\Delta(x_0)} \frac{r(t,y_h)(x_h-y_h)}{|x_h-y_h|^3}\ \dy_h\Big| <\delta.
\]
Whence (\ref{G4}).

Therefore we proved that
\[
|K_3 +R_9| \leq  \delta \| \vc U-\vu_\ep\|_{L^6(\Omega;\RRR^3)} + C(\delta; r, \vc{V},\Theta) \int_\Omega {\mathcal E}(\vr_\ep,\vt_\ep|r,\Theta)\, \dx + H_\ep,
\]
where $H_\ep = o(\ep)$ as $\ep \to 0^+$.

Plugging all the previous estimates into (\ref{REI}) we get
\[
\int_{\Omega}
\left(
\frac{1}{2}\ \vr_{\varepsilon} |\vc u_{\varepsilon}-\vc U|^2
+{\mathcal E}\left(\vr_{\varepsilon},\vt_{\varepsilon}|r ,\Theta\right)\right)(\tau,\cdot)\ \dx
\]
\[
+\int_0^{\tau}\int_{\Omega}
\Big(
\frac{\Theta}{\vt_{\varepsilon}}
 \tn{S}(\vt_{\varepsilon},\Grade\vu_{\varepsilon}) : \Grade \vu_{\varepsilon}
-\tn{S}(\Theta,\Grade\vc U) : (\Grade \vu_{\varepsilon}-\Grade \vc U)
-\tn{S}(\vt_{\varepsilon},\Grade\vu_{\varepsilon}) : \Grade \vc U
\Big)\ \dx\ \dt
\]
\[
 +\int_0^{\tau}\int_{\Omega}
\Big(
\frac{\vc{q}(\vt_{\varepsilon},\Grade\vt_{\varepsilon}) \cdot \Grade \Theta }{\vt_{\varepsilon}}
-\frac{\Theta}{\vt_{\varepsilon}}\ \frac{\vc{q}(\vt_{\varepsilon},\Grade\vt_{\varepsilon}) \cdot \Grade \vt_{\varepsilon} }{\vt_{\varepsilon}}
\Big)
\ \dx\ \dt
\]
\[
\leq
\int_{\Omega}
\frac{1}{2}
\left(
 \vr_{0,\varepsilon} |\vc u_{0,\varepsilon}- \vc U(0,\cdot)|^2
+{\mathcal E}
\left(
\vr_{0,\varepsilon},\vt_{0,\varepsilon}|r(0,\cdot),\Theta(0,\cdot)
\right)
\right)\ \dx + H_\ep
\]
\[
+
\int_0^{\tau}
\Big[
\delta\|\vc V-\vc u_{\varepsilon}\|^2_{W^{1,2}(\Omega;\RRR^3)}
+{\mathcal C}(\delta;r,\vc V,\Theta)
\int_{\Omega}
\Big(
\frac{1}{2}\ \vr_{\varepsilon} |\vc u_{\varepsilon}-\vc U|^2
+{\mathcal E}\left(\vr_{\varepsilon},\vt_{\varepsilon}|r ,\Theta\right)
\Big)
\ \dx
\Big]\ \dt
\]
\[
+
\int_0^\tau\int_{\Omega}
\Big(  p(r ,\Theta)-p(\vr_{\varepsilon},\vt_{\varepsilon})\Big)\Divh \vc V\ \dx\ \dt \]
\[
+\int_0^\tau\int_{\Omega}
\Big(
\left(
1-\frac{ \vr_{\varepsilon}}{r}
\right)
\big(
\partial_t p(r ,\Theta)
+ \vc V\cdot \Gradh  p(r ,\Theta)
\Big)
\ \dx \ \dt
\]
\[
-\int_0^\tau\int_{\Omega}
  r \Big(
\partial_{\vr}s(r ,\Theta)(\vr_{\varepsilon} -r)
+\partial_{\vt}s(r ,\Theta)(\vt_{\varepsilon} -\Theta)
\Big)
\big(
\partial_t \Theta+\vc V\cdot  \Gradh \Theta\big) \ \dx \ \dt,
\]
where $H_\ep $ is as above.
We denote
\[
{\mathcal A}:=\int_{\Omega}
\big(  p(r ,\Theta)-p(\vr_{\varepsilon},\vt_{\varepsilon})\big)\Divh \vc V
\ \dx
+\int_{\Omega}
\left(
1-\frac{ \vr_{\varepsilon}}{r}
\right)
\big(
\partial_t p(r ,\Theta)
+ \vc V\cdot \Gradh  p(r ,\Theta)
\big)
\ \dx
\]
\[
-\int_{\Omega}
  r \Big(
\partial_{\vr}s(r ,\Theta)(\vr_{\varepsilon} -r)
+\partial_{\vt}s(r ,\Theta)(\vt_{\varepsilon} -\Theta)
\Big)
\big(
\partial_t \Theta+\vc V\cdot  \Gradh \Theta
\big)
 \ \dx
\]
\[
=\int_{\Omega}
\big(  p(r ,\Theta)-p(\vr_{\varepsilon},\vt_{\varepsilon})\big)\Divh \vc V
\ \dx
+\int_{\Omega}
r
\left(
\Theta- \vt_{\varepsilon}
\right)
\partial_{\vt}s(r ,\Theta)
\big(
\partial_t \Theta
+ \vc V\cdot \Gradh  \Theta
\big)
\ \dx
\]
\[
-\int_{\Omega}
   \big(r-\vr_{\varepsilon}\big)
\partial_{\vr}p(r ,\Theta)\Divh \vc V\ \dx,
\]
where we used the Gibbs (\ref{m4a}) and Maxwell (\ref{i5}) relations and the continuity equation (\ref{t1}).
The second term in the right-hand side rewrites as follows
\[
\int_{\Omega}
r
\left(
\Theta- \vt_{\varepsilon}
\right)
\partial_{\vt}s(r ,\Theta)
\big(
\partial_t \Theta
+ \vc V\cdot \Gradh  \Theta
\big)
\ \dx
\]
\[
=
\int_{\Omega}
r
\left(
\Theta- \vt_{\varepsilon}
\right)
\big(
\partial_t s(r ,\Theta)
+ \vc V\cdot \Gradh  s(r ,\Theta)
\big)
\ \dx
-\int_{\Omega}
   (\Theta-\vt_{\varepsilon})
\partial_{\vt}p(r ,\Theta)\Divh \vc V\ \dx
\]
\[
=
\int_{\Omega}
\left(
\Theta- \vt_{\varepsilon}
\right)
\Big[
\frac{1}{\Theta}
\Big(
 \tn{S}_h(\Theta,\Gradh\vc V) : \Gradh \vc V
-\frac{\vc{q}_h(\Theta,\Gradh \Theta) \cdot \Gradh \Theta }{\Theta}
\Big)
-\Divh
\Big(\frac{\vc{q}_h(\Theta,\Gradh \Theta) }{\Theta}
\Big)
\Big]\ \dx
\]
\[
-\int_{\Omega}
   (\Theta-\vt_{\varepsilon})
\partial_{\vt}p(r ,\Theta)\Divh \vc V\ \dx,
\]
where we used the same identities as above, together with the entropy balance (\ref{t3}). Hence
\[
{\mathcal A}=\int_{\Omega}
\Big(
 p(\vr,\Theta)
-p(\vr_{\varepsilon},\vt_{\varepsilon})
+\partial_{\vr} p(r,\Theta)(\vr_{\varepsilon} -r)
+\partial_{\vt} p(r,\Theta)(\vt_{\varepsilon} -\Theta)
\Big)
\Divh \vc V\ \dx
\]
\[
+\int_{\Omega}
\left(
\Theta- \vt_{\varepsilon}
\right)
\Big[
\frac{1}{\Theta}
\Big(
 \tn{S}_h(\Theta,\Gradh\vc V) : \Gradh \vc V
-\frac{\vc{q}_h(\Theta,\Gradh \Theta) \cdot \Grad \Theta }{\Theta}
\Big)
-\Divh
\Big(\frac{\vc{q}_h(\Theta,\Grad \Theta) }{\Theta}
\Big)
\Big]\ \dx.
\]
Observing that
\[
\left|\int_{\Omega}
\Big(
 p(r ,\Theta)
-p(\vr_{\varepsilon},\vt_{\varepsilon})
+\partial_{\vr} p(r ,\Theta)(\vr_{\varepsilon} -r)
+\partial_{\theta} p(r ,\Theta)(\vt_{\varepsilon} -\Theta)
\Big)
\Divh \vc V\ \dx
 \right|
\]
\[
\leq
C\|\Divh\vc V\|_{L^{\infty}(\Omega)}
\int_{\Omega}{\mathcal E}\left(\vr_{\varepsilon},\vt_{\varepsilon}|r ,\Theta\right)\ \dx,
\]
we see that (\ref{REI}) reduces finally to
\[
\int_{\Omega}
\left(
\frac{1}{2}\ \vr_{\varepsilon} |\vc u_{\varepsilon}-\vc U|^2
+{\mathcal E}\left(\vr_{\varepsilon},\vt_{\varepsilon}|r ,\Theta\right)
\right)(\tau,\cdot)\ \dx
\]
\[
+\int_0^{\tau}\int_{\Omega}
\Big(
\frac{\Theta}{\vt_{\varepsilon}}
 \tn{S}(\vt_{\varepsilon},\Grade\vu_{\varepsilon}) : \Grade \vu_{\varepsilon}
-\tn{S}(\Theta,\Grade \vc U) : (\Grade \vu_{\varepsilon}-\Grade \vc U)
-\tn{S}(\vt_{\varepsilon},\Grade\vu_{\varepsilon}) : \Grade \vc U
\Big)\ \dx\ \dt
\]
\[+\int_0^{\tau}\int_{\Omega}
\frac{\Theta-\vt_\ep}{\Theta}\tn{S}_h(\Theta,\Gradh\vc V) : \Gradh \vc V
\ \dx\ \dt
\]\[
 +\int_0^{\tau}\int_{\Omega}
\Big(
\frac{\vc{q}(\vt_{\varepsilon},\Grade\vt_{\varepsilon}) \cdot \Grade \Theta }{\vt_{\varepsilon}}
-\frac{\Theta}{\vt_{\varepsilon}}\ \frac{\vc{q}(\vt_{\varepsilon},\Grade \vt_{\varepsilon}) \cdot \Grade \vt_{\varepsilon} }{\vt_{\varepsilon}}
\Big)
\ \dx\ \dt
\]
\[
 +\int_0^{\tau}\int_{\Omega}
\Big(
(\Theta-\vt_{\varepsilon})\frac{\vc{q}_h(\Theta,\Gradh\Theta) \cdot \Gradh \Theta }{\Theta^2}
+\frac{\vc{q}(\Theta,\Grade\Theta) \cdot \Grade (\vt_{\varepsilon}-\Theta) }{\Theta}
\Big)
\ \dx\ \dt
\]
\[
\leq
\int_{\Omega}
\frac{1}{2}
\left(
 \vr_{0,\varepsilon} |\vc u_{0,\varepsilon}- \vc U(0,\cdot)|^2
+{\mathcal E}
\left(
\vr_{0,\varepsilon},\vt_{0,\varepsilon}|r(0,\cdot),\Theta(0,\cdot)
\right)
\right)\ \dx + H_\ep
\]
\bFormula{EE2}
+\int_0^{\tau}
\Big[
\delta\|\vc U-\vc u_{\varepsilon}\|^2_{W^{1,2}(\Omega;\RRR^3)}
+C(\delta;r,\vc V,\Theta)
\int_{\Omega}{\mathcal E}\left(\vr_{\varepsilon},\vt_{\varepsilon}|r ,\Theta\right)\ \dx
\Big]\ \dt.
\eF
In order to finish the proof, it remains to study carefully the terms in the left-hand side.

\subsection{Dissipative terms}
We will consider more carefully only the terms with the stress tensor, referring to \cite{FN} for more details, especially in the case of the terms with the heat flux.
In the spirit of \cite{FN}, using the special structure of the stress tensor (see (\ref{m7})),
we can write
\begin{equation}
 \tn{S}(\vt_\ep,\Grade \vu_\ep) = \tn{S}^0(\Grade \vu_\ep)+\tn{S}^1(\vt_\ep,\Grade \vu_\ep),
 \end{equation}
 where
 \[\tn{S}^0(\Grade\vu_\ep)= \mu_0\left(\Grade \vu_\ep+ (\Grade \vu_\ep)^T - \frac 23 \Dive \vu_\ep \, \tn{I}\right),\]
 \[\tn{S}^1(\vt_\ep,\Grade \vu_\ep)=\mu_1 \vt_\ep \left(\Grade \vu_\ep+ (\Grade \vu_\ep)^T - \frac 23 \Dive \vu_\ep \, \tn{I}\right).\]
Then
 \[
\frac{\Theta}{\vt_{\varepsilon}}
 \tn{S}^1(\vt_{\varepsilon},\Grade\vu_{\varepsilon}) : \Grade \vu_{\varepsilon}
-\tn{S}^1(\Theta,\Grade\vc U) : (\Grade \vu_{\varepsilon}-\Grade \vc U)
-\tn{S}^1(\vt_{\varepsilon},\Grade\vu_{\varepsilon}) : \Grade \vc U
\]
\[+
\Big(\frac{\vt_{\varepsilon}- \Theta}{\Theta}\tn{S}_h^1(\Theta,\Gradh\vc V) : \Gradh \vc V
\Big) \]
\[= \Theta \Big (\frac{\tn{S}^1(\vt_{\varepsilon},\Grade\vu_{\varepsilon})}{\vt_{\varepsilon}} -\frac{\tn{S}^1(\Theta,\Grade\vc U)}{\Theta}\Big):  (\Grade \vu_{\varepsilon}-\Grade \vc U)\]
\[+(\Theta - \vt_{\varepsilon})\Big (\frac{\tn{S}^1(\vt_{\varepsilon},\Grade\vu_{\varepsilon})}{\vt_{\varepsilon}} -\frac{\tn{S}^1(\Theta,\Grade\vc U)}{\Theta}\Big):\Grade	 \vc U.\]

In the first term we use the Korn inequality and we split the second term to the essential and the residual parts. As above, we get bounds
\[\Big|\int_\Omega (\Theta - \vt_{\varepsilon})\Big (\frac{\tn{S}^1(\vt_{\varepsilon},\Grade\vu_{\varepsilon})}{\vt_{\varepsilon}} -\frac{\tn{S}^1(\Theta,\Grade\vc U)}{\Theta}\Big):\Grade	 \vc U \ \dx \Big|\]
\[\leq \delta \|\vu_\ep  -\vc U\|^2_{W^{1,2}(\Omega;\RRR^3)}+C(\delta)\int _{\Omega}{\mathcal E}\left(\vr_{\varepsilon},\vt_{\varepsilon}|r ,\Theta\right)\ \dx.\]

The part of the stress tensor containing $\tn{S}^0$  can be rewritten in the following way. First, for $0<\Theta \leq \vt_\ep$
\[
\frac{\Theta}{\vt_{\varepsilon}}
 \tn{S}^0(\Grade\vu_{\varepsilon}) : \Grade \vu_{\varepsilon}
-\tn{S}^0(\Grade\vc U) : (\Grade \vu_{\varepsilon}-\Grade \vc U)
-\tn{S}^0(\Grade\vu_{\varepsilon}) : \Grade \vc U
\]
\[+\frac{\vt_{\varepsilon}- \Theta}{\Theta}\tn{S}^0_h(\Gradh\vc V) : \Gradh \vc V
  \]
\[\geq \frac{\Theta}{\vt_{\varepsilon}} \big(\tn{S}^0(\Grade\vu_{\varepsilon})-\tn{S}^0(\Grade \vc U)\big):\Grade (\vu_{\varepsilon}- \vc U)+ \Theta\big ( \frac{1}{\vt_{\varepsilon}}-\frac{1}{\Theta}\big ) \tn{S}^0(\Grade \vc U):\Grade (\vu_{\varepsilon}- \vc U)\]
\[+ \frac{\vt_\ep-\Theta}{\vt_\ep}\big(\tn{S}^0(\Grade\vc U)-\tn{S}^0(\Grade\vu_\ep)\big):\Grade\vc{U}. \]
As $\frac{1}{\vt_\ep}\leq \frac{1}{\Theta}$, the last two terms can be, after integrating over $\Omega$, estimated by
\[\delta \|\vu_\ep-\vc{U}\|_{W^{1,2}(\Omega;\RRR^3)}^2 + C(\delta) \int_\Omega {\mathcal E}(\vr_\ep,\vt_\ep|r,\Theta)\ \dx.
\]
Next, for $0 < \vt_\ep\leq \Theta$,
\[
\frac{\Theta}{\vt_{\varepsilon}}
 \tn{S}^0(\Grade\vu_{\varepsilon}) : \Grade \vu_{\varepsilon}
-\tn{S}^0(\Grade\vc U) : (\Grade \vu_{\varepsilon}-\Grade \vc U)
-\tn{S}^0(\Grade\vu_{\varepsilon}) : \Grade \vc U
\]
\[+\frac{\vt_{\varepsilon}- \Theta}{\Theta}\tn{S}^0(\Gradh\vc V) : \Gradh \vc V
\]
\[\geq \big(\tn{S}^0(\Grade\vu_\ep)-\tn{S}^0(\Grade\vc U)\big):\Grade(\vu_\ep-\vc{U}) + \frac{\Theta-\vt_\ep}{\Theta}\big(\tn{S}^0(\Grade \vu_\ep):\Grade\vu_\ep - \tn{S}^0_h (\Gradh \vc V): \Gradh \vc{V} \big).
\]
As $\Grade \vu_{\ep} \mapsto \tn{S}^0(\Grade \vu_\ep):\Grade\vu_\ep$ is convex, we have
\[
\frac{\Theta-\vt_\ep}{\Theta}\big(\tn{S}^0(\Grade \vu_\ep):\Grade\vu_\ep - \tn{S}^0_h (\Gradh \vc V): \Gradh \vc{V}\big) \geq \frac{\Theta-\vt_\ep}{\Theta}\tn{S}^0(\Grade\vc{U}):\Grade(\vu_\ep-\vc U).
\]
After integration over $\Omega$, this term can be controlled on the right-hand side by
\[\delta \|\vu_\ep-\vc{U}\|_{W^{1,2}(\Omega;\RRR^3)}^2 + C(\delta) \int_\Omega {\mathcal E}(\vr_\ep,\vt_\ep|r,\Theta)\ \dx.
\]
Therefore, summing up,
\[
\int_{\Omega}
\left(
\frac{1}{2}\ \vr_{\varepsilon} |\vc u_{\varepsilon}-\vc U|^2
+{\mathcal E}\left(\vr_{\varepsilon},\vt_{\varepsilon}|r ,\Theta\right)
\right)(\tau,\cdot)\ \dx
+k_1\int_0^{\tau}\int_{\Omega}
\big|
\Grade\vu_{\varepsilon} - \Grade \vc U
\big|^2\ \dx\ \dt\]
\[
 +\int_0^{\tau}\int_{\Omega}
\Big(
\frac{\vc{q}(\vt_{\varepsilon},\Grade\vt_{\varepsilon}) \cdot \Grade \Theta }{\vt_{\varepsilon}}
-\frac{\Theta}{\vt_{\varepsilon}}\ \frac{\vc{q}(\vt_{\varepsilon},\Grade \vt_{\varepsilon}) \cdot \Grad \vt_{\varepsilon} }{\vt_{\varepsilon}}
\Big)
\ \dx\ \dt
\]
\[
 +\int_0^{\tau}\int_{\Omega}
\Big(
(\Theta-\vt_{\varepsilon})\frac{\vc{q}_h(\Theta,\Gradh\Theta) \cdot \Gradh \Theta }{\Theta^2}
+\frac{\vc{q}_h(\Theta,\Grade\Theta) \cdot \Grade (\vt_{\varepsilon}-\Theta) }{\Theta}
\Big)
\ \dx\ \dt
\]
\[
\leq
\int_{\Omega}
\frac{1}{2}
\left(
 \vr_{0,\varepsilon} |\vc u_{0,\varepsilon}- \vc U(0,\cdot)|^2
+{\mathcal E}
\left(
\vr_{0,\varepsilon},\vt_{0,\varepsilon}|r(0,\cdot),\Theta(0,\cdot)
\right)
\right)\ \dx
\]
\bFormula{EE2a}
+\int_0^{\tau}
\Big[
\delta\|\vc U-\vc u_{\varepsilon}\|^2_{W^{1,2}(\Omega;\RRR^3)}
+C(\delta;r,\vc V,\Theta)
\int_{\Omega}{\mathcal E}\left(\vr_{\varepsilon},\vt_{\varepsilon}|r ,\Theta\right)\ \dx
\Big]\ \dt
.
\eF

For the terms connected with the heat conductivity, we proceed exactly as in \cite{FN}. We arrive at the following inequality
\[
\int_{\Omega}
\left(
\frac{1}{2}\ \vr_{\varepsilon} |\vc u_{\varepsilon}-\vc U|^2
+{\mathcal E}\left(\vr_{\varepsilon},\vt_{\varepsilon}|r ,\Theta\right)
\right)(\tau,\cdot)\ \dx
+k_1\int_0^{\tau}\int_{\Omega}
\big|
\Grade\vu_{\varepsilon} - \Grade\vc U
\big|^2\ \dx\ \dt
\]
\[
 +k_2\int_0^{\tau}\int_{\Omega}
\big|
\Grade\vt_{\varepsilon} - \Grade \Theta
\big|^2\ \dx\ \dt
 +k_3\int_0^{\tau}\int_{\Omega}
\big|
\Grade\log\vt_{\varepsilon} - \Grade\log \Theta
\big|^2\ \dx\ \dt
\]
\[
\leq
\int_{\Omega}
\Big(
\frac{1}{2}
\left(
 \vr_{0,\varepsilon} |\vc u_{0,\varepsilon}- \vc U(0,\cdot)|^2
+{\mathcal E}
\left(
\vr_{0,\varepsilon},\vt_{0,\varepsilon}|r(0,\cdot),\Theta(0,\cdot)
\right)
\right)
\ \dx + H_\ep
\]
\bFormula{EE3}
+k_4\int_0^{\tau}\int_{\Omega}
\Big(
\frac{1}{2}\ \vr_{\varepsilon}\big|\vc u_{\varepsilon}-\vc U\big|^2
+{\mathcal E}\left(\vr_{\varepsilon},\vt_{\varepsilon}|r ,\Theta\right)
\Big)\ \dx\ \dt,
\eF
where the positive constants $k_j$ depend on $(r,\vc V,\Theta)$ through the norms
 involved in Theorem \ref{main}, and $H_\ep \to 0$ as $\ep \to 0^+$. This inequality, after application of the Gronwall lemma, finishes the proof of Theorem \ref{main}.

\smallskip

{\bf Acknowledgement:} Bernard Ducomet is partially supported by the ANR project INFAMIE (ANR-15-CE40-0011).
\v S\' arka Ne\v casov\' a and Matteo Caggio acknowledge the support of 
the GA\v CR (Czech Science Foundation) project No. 16-03230S in the framework of RVO: 67985840. Part of article was written during  stay of \v S\' arka Ne\v casov\' a at CEA. She would like to thank to Prof. Ducomet for his hospitality. 
 The work of Milan Pokorn\'y was supported by the GA\v CR (Czech Science Foundation) project No. 16-03230S.

\end{document}